\documentclass{amsart}
\usepackage{tabu}
\usepackage[margin=95pt]{geometry}
\usepackage[round,authoryear]{natbib}
\usepackage[hidelinks]{hyperref}
\usepackage{enumitem}
\allowdisplaybreaks
\usepackage{mathrsfs}
\usepackage{comment}
\usepackage{MnSymbol}
\usepackage{accents}

\DeclareMathOperator{\quo}{Q}
\DeclareMathOperator{\rem}{R}
\DeclareMathOperator{\lcm}{lcm}
\renewcommand{\pmb}[1]{\underaccent{\tilde}{#1}}

\newcommand{\Z}{\mathbb{Z}}
\newcommand{\R}{\mathbb{R}}
\newcommand{\cK}{\mathscr{K}}
\newcommand{\cI}{\mathscr{I}}
\newcommand{\cW}{\mathcal{W}}
\newcommand{\cU}{\mathcal{U}}
\DeclareMathOperator{\fJ}{Bad}
\newcommand{\fS}{\mathfrak{S}}

\begin {document}

\newtheorem{thm}{Theorem}
\newtheorem{alg}{Algorithm}
\newtheorem{lem}{Lemma}[section]
\newtheorem{prop}[lem]{Proposition}

\title[sums of seven cubes]{Every integer greater than $454$ is\\
 the sum
of at most seven positive cubes}

\author{Samir Siksek}

\address{Mathematics Institute\\
	University of Warwick\\
	United Kingdom}


\date{\today}
\thanks{
The 
author is supported by an EPSRC Leadership Fellowship EP/G007268/1,
and EPSRC {\em LMF: L-Functions and Modular Forms} Programme Grant
EP/K034383/1.
}


\begin{abstract}
A long-standing conjecture states that
every positive integer 
other than
\begin{gather*}
15,\; 22,\; 23,\; 50,\; 114,\; 167,\; 175,\; 186,\; 212,\\
231,\; 238,\; 239,\; 303,\; 364,\; 420,\; 428,\; 454
\end{gather*}
is a sum of at most
seven positive cubes. 
This was first 
observed by Jacobi in 1851 on the basis of
extensive calculations performed by the famous computationalist Zacharias Dase.
We complete the proof of this conjecture,
building on previous work of Linnik, Watson, 
McCurley, Ramar\'{e},
Boklan, Elkies, and many others.
\end{abstract}
\maketitle


\section{Historical Introduction}
In 1770, Edward Waring stated in his Meditationes Algebraic\ae, 
\begin{quote}
\emph{Omnis integer numerus vel est cubus, vel e duobus, tribus, 4, 5, 6, 7,
8, vel novem cubis compositus,  \dots}
\end{quote}
Waring's assertion, can be concisely reformulated as: 
every positive integer is the sum of 
nine non-negative cubes. Henceforth, \textbf{by a cube we shall mean
a non-negative cube}.
In the 19th century, numerical experimentation led to 
refinements of Waring's assertion for sums of cubes. 
As noted by \citet{DicksonExtension},
\emph{\lq\lq At the request of Jacobi, the famous computer Dase constructed a table
showing the least number of positive cubes whose sum is any $p<12000$\rq\rq}. 
In an influential Crelle paper, 
\citet{Jacobi} made a series of observations 
based on Dase's table: every positive integer 
other than $23$ and $239$ is the sum of
eight cubes, every integer $>454$ is the sum
of seven cubes, and every integer $>8042$
is the sum of six cubes. Jacobi believed that every sufficiently large 
integer is the sum of five cubes, whilst recognizing that the 
cut-off point must be far beyond Dase's table,
and he wondered if the same is true for sums of four cubes. He noted
that integers $\equiv 4$, $5 \pmod{9}$ cannot be sums of three cubes.
Later computations by \citet{Romani} convincingly suggest
that every integer $>1\,290\,740$ is the sum of five cubes,
and by 
\citet*{Desh}
that every integer $>7\,373\,170\,279\,850$ is the sum of four cubes.

\bigskip

Progress towards proving these observations of Waring, Jacobi and others 
has been exceedingly slow. \citet{Maillet} showed that twenty-one 
cubes are enough to represent
every positive integer. At the heart of Maillet's proof is an
 idea crucial to 
virtually all
future developments; the identity $(r+x)^3+(r-x)^3=2 r^3+6 r x^2$ allows one to
reformulate the problem of representing an integer as the sum of a (certain number of) cubes in terms of 
representing a related integer as the sum of (a smaller number of) squares.
Exploiting this idea, \citet{Wieferich} proved Waring's assertion
(Wieferich's proof
had a mistake that was corrected by \citet{Kempner}). 
In fact, the theoretical part of Wieferich's proof showed that all integers 
exceeding $2.25 \times 10^9$ are sums of nine cubes. 
Completing the proof required appealing to a table of \citet{Sterneck} 
(who extended Dase's table to $40\,000$), and applying what is now known 
as the greedy algorithm to reach the bound.

\bigskip

Soon thereafter, \citet{Landau}
showed that 
every sufficiently large integer is the sum of eight cubes.
This was made effective by \citet{Baer}, who showed that every 
integer $\ge 14.1 \times 233^6 \approx 2.26 \times 10^{15}$ is the sum of eight cubes. 
\citet{DicksonEight} completed the proof
of Jacobi's observation that all positive 
integers other than $23$ and $239$
are sums of eight cubes.
Remarkably, Dickson's proof relied on extending von Sterneck's table to $123\,000$ (with the
help of his assistant, Miss Evelyn Garbe) and then 
applying the greedy algorithm to
reach Baer's bound.

\bigskip

In \citeyear{Linnik} Linnik showed that every sufficiently large
integer is the sum of seven cubes. A 
substantially simpler proof (though still ineffective)  
was given by \citet{Watson}. Linnik's seven cubes theorem was 
first
made effective by \citet{McCurley},
who showed that it is true for integers $>\exp(\exp(13.94))$.
Ramar\'{e} improved this to $\exp(205000)$ 
in \citeyear{RamareBig}
and finally 
to $\exp(524)\approx 3.72 \times 10^{227}$
in \citeyear{RamareSmall}. 
This bound is way beyond computer searches combined with
the greedy algorithm. 
In \citet{Desh}, it is shown that every integer between $1\,290\,741$ and $10^{16}$
is a sum of five cubes. As observed by \citet{RamareSmall},
combining this with
the greedy algorithm  \citep[Lemma 3]{BRZ}, we can easily deduce 
that every 
integer $455 \le N \le \exp(78.7) \approx 1.51 \times 10^{34}$ is the sum of seven cubes.

\bigskip

There has been a number of partial results concerning sums of seven cubes.
\citet{BRZ} show that
every non-negative integer which is a cubic residue modulo $9$ and an invertible
cubic residue modulo $37$ is a sum of $7$ cubes.
\citet{BE} show that every multiple of $4$ greater than $454$ 
is the sum of seven cubes, whilst \citet{Elkies} shows
the same for integers $\equiv 2 \pmod{4}$.

\medskip

In this paper we complete the proof of Jacobi's seven cubes
conjecture, building on the aforementioned great works.
\begin{thm}\label{thm:main}
Every positive integer 
other than
\[
15,\; 22,\; 23,\; 50,\; 114,\; 167,\; 175,\; 186,\; 212,\; 
231,\; 238,\; 239,\; 303,\; 364,\; 420,\; 428,\; 454
\]
is the sum of seven cubes.
\end{thm}

The programs that accompany this paper are available from
\url{http://tinyurl.com/zlaeweo}.
It is a pleasure to thank 
Alex Bartel, Tim Browning, John Cremona, Roger Heath-Brown and Trevor Wooley
for stimulating discussions.

\section{The Main Criterion}\label{sec:criterion}
Let $\cK=\exp(524)$ and $\cK^\prime=\exp(78.7)$.
By the results of \citet{RamareSmall}
and of \citet{Desh}, it is sufficient to prove that every
integer $\cK^\prime \le N \le \cK$ is the sum of 
seven cubes. 
The results of \cite{BE} and \cite{Elkies} allow us
to restrict ourselves to odd integers $N$
(our method can certainly be adapted to deal with 
even integers, but
restricting ourselves to odd integers brings 
coherence to our exposition).
In this section we give a criterion (Proposition~\ref{prop:sevenK}) 
for  all odd integers $N$
in a range $K_1 \le N \le K_2$
to be sums of seven  cubes.
Most of the remainder of the paper is devoted to
showing that this criterion holds for each of
the ranges 
$(9/10)^{n+1} \cK \le N \le (9/10)^n \cK$ with
$0 \le n \le 4226$. This will complete
the proof of Theorem~\ref{thm:main}
as $(9/10)^{4227} \cK \approx 1.42 \times 10^{34}$
and $\cK^\prime \approx 1.51 \times 10^{34}$.

\begin{thm}[Gauss, Legendre]\label{lem:threesq} 
Let $k \ge 0$ be an even integer. There exist 
integers
$x$, $y$, $z$ such that
\begin{equation}\label{eqn:threesq}
x^2+x+y^2+y+z^2+z=k.
\end{equation}
\end{thm}
\begin{proof}
Dividing by $2$ we see that this is in fact the famous theorem, due to Gauss,
that every non-negative integer is the sum of three triangular numbers.
Alternatively, we can rewrite \eqref{eqn:threesq} as
\begin{equation}\label{eqn:threesq2}
(2x+1)^2+(2y+1)^2+(2z+1)^2=4k+3.
\end{equation}
As $k$ is even, $4k+3 \equiv 3 \pmod{8}$;
by a theorem of Legendre,
 every positive integer $\equiv 3 \pmod{8}$
is the sum of three odd squares.
\end{proof}

Throughout this section
$m$ will denote a positive integer satisfying the conditions
\begin{enumerate}
\item[(i)] $m$ is a squarefree,
\item[(ii)] $3 \mid m$,
\item[(iii)] every prime divisor of $m/3$ is $\equiv 5 \pmod{6}$.
\end{enumerate}
Observe that $m \equiv 3 \pmod{6}$. Moreover, for any 
integer $N$, there is a unique integer $t \in [0,m)$ such 
that $N \equiv 8 t^3 \pmod{m}$. 
Our starting point is a modified version of Lemma~3 of \cite{Watson}.
\begin{lem}\label{lem:sevenK}
Let $0<K_1 < K_2$ be real numbers. Let $m$ be a positive integer  
satisfying (i)--(iii) above.
Let
$\varepsilon_m$, $\delta_m$ be real numbers satisfying
\begin{enumerate}
\addtolength\itemsep{1mm}
\item[(iv)] $0 \le \varepsilon_m < \delta_m \le 1$,
\item[(v)] $K_1 \ge (8 \delta_m^3+1/36) m^3+3m/4$,
\item[(vi)] $K_2 \le (8 \varepsilon_m^3+1/18) m^3+ m/2$.
\end{enumerate}
Let $K_1 \le N \le K_2$ be an odd integer. Suppose 
$N \equiv 8 t^3 \pmod{m}$ with 
$t \in [\varepsilon_m \cdot m, \delta_m \cdot m)$.
Then $N$ is the sum of seven non-negative cubes.
\end{lem}
\begin{proof}
Write $m=6r+3$.
Let 
\begin{equation}\label{eqn:k}
k=\frac{N-8 t^3}{m}-(r^2+r+1).
\end{equation}
The quantity $k$ is an integer as $N \equiv 8 t^3 \pmod{m}$, 
and even as $(N-8t^3)/m$ and $r^2+r+1$ are both odd.
Observe that
\begin{align*}
k & > \frac{N-8\delta_m^3 \cdot m^3}{m}-(r^2+r+1) \qquad \text{as $t < \delta_m \cdot m$}\\
& \ge \frac{K_1-8\delta_m^3 \cdot m^3}{m}-(r^2+r+1) \qquad \text{as $N \ge K_1$}\\
& = \frac{K_1-8\delta_m^3 \cdot  m^3}{m}-\frac{m^2}{36}-\frac{3}{4} 
\qquad \text{substituting $r=(m-3)/6$}\\
  & \ge 0 \qquad \text{by (v)}.
\end{align*}
As $k$ is non-negative and even, by the Gauss--Legendre theorem,
there exist integers $x$, $y$, $z$ 
satisfying \eqref{eqn:threesq}.
We shall make use of the identity
\begin{multline}\label{eqn:identity}
(r+1+x)^3+(r-x)^3+(r+1+y)^3+(r-y)^3+(r+1+z)^3+(r-z)^3
=\\
(6r+3) (r^2+r+1+x^2+x+y^2+y+z^2+z).
\end{multline}
From the definition of $k$ in \eqref{eqn:k} and the fact
that $m=6r+3$, we see that $N-8t^3$ is equal to the right-hand side
of the identity \eqref{eqn:identity}. Hence
\[
N=
(r+1+x)^3+(r-x)^3+(r+1+y)^3+(r-y)^3+(r+1+z)^3+(r-z)^3+
(2 t)^3.
\]
To complete the proof it is enough to show that
these cubes are non-negative, or equivalently that
\[
-r-1
\le
x,\; y,\; z \le r.
\]
This is equivalent to showing that
\[
-(2r+1) \le 2x+1, \; 2y+1, \; 2z+1 \le 2r+1.
\]
Now $(2y+1)^2$, $(2z+1)^2 \ge 1$ and so from
\eqref{eqn:threesq2}, we have $(2x+1)^2 \le 4k+1$.
 It is therefore enough to show that
$4k+1 \le (2r+1)^2$ or equivalently that
$k \le r^2+r$. The following inequalities complete the proof:
\begin{align*}
k -r^2-r& 
=\frac{N-8 t^3}{m}-(2r^2+2r+1) \qquad \text{from \eqref{eqn:k}}\\
& \le \frac{N-8\varepsilon_m^3 \cdot m^3}{m}-(2r^2+2r+1) \qquad \text{as $t \ge \varepsilon_m \cdot m$}\\
& \le \frac{K_2-8\varepsilon_m^3 \cdot m^3}{m}-(2r^2+2r+1) \qquad \text{as $N \le K_2$}\\
& \le \frac{K_2-8\varepsilon_m^3 \cdot m^3}{m}-\frac{m^2}{18}-\frac{1}{2} \qquad 
  \text{substituting $r=(m-3)/6$}\\
& \le 0 \qquad \text{by (vi)} \, .
\end{align*}
\end{proof}
This simple-minded lemma has one serious flaw. The condition $K_1<K_2$
together with conditions (iv), (v), (vi), imply 
\[
\delta_m^3 < \varepsilon_m^3+1/288-1/(32 m^2) \, .
\] 
In particular, this forces the interval $[\varepsilon \cdot m, \delta \cdot m)$
to have length $< m/\sqrt[3]{288} \approx 0.15m$. On the other hand,
the integer $t$ appearing in the lemma (which is the cube-root of $N/8$
modulo $m$) can be any integer in the interval $[0,m)$. Thus the lemma
only treats a small fraction of the odd integers $K_1 \le N \le K_2$.
Our key innovation over the works mentioned in the introduction
is to use not just one value of $m$, but many of them simultaneously.
Each value of $m$ will give some information about those odd
integers $K_1 \le N \le K_2$ that cannot be expressed as sums
of seven cubes; collecting this information will allow us to deduce
a contradiction.

\medskip

Let $x$ be a real number and $m$ be a positive integer.
Define the \textbf{quotient} and \textbf{remainder} 
obtained on dividing $x$ by $m$ as
\[
 \quo(x,m)=\lfloor x/m\rfloor, \qquad \rem(x,m) =x-\quo(x,m) \cdot m.
\]
In particular, $\rem(x,m)$ belongs to the half-open interval $[0,m)$.
If $x\in \Z$ then $\rem(x,m)$ is the usual remainder on dividing by $m$,
and $x \equiv \rem(x,m) \pmod{m}$.
Let $\varepsilon$ and $\delta$ be real numbers satisfying $0 \le \varepsilon < \delta \le 1$. 
Define
\begin{equation}
\label{eqn:fJ}
\fJ(m,\varepsilon,\delta) 
=\Big\{
\;  x\in \R  \; : \; \rem(x,m) 
\in {[0,m)} \setminus {[\varepsilon \cdot m,\delta \cdot m)} 
\; \Big\} 
=\bigcup_{k=-\infty}^{\infty} km+ \left( {[0,m)} \setminus {[\varepsilon \cdot m,\delta \cdot m)}\right).
\end{equation}
The reader will observe, in Lemma~\ref{lem:sevenK}, if $N$ is not
the sum of seven cubes, then $t \in \fJ(m,\varepsilon,\delta)$,
which explains our choice of the epiphet \lq bad\rq.
Given a set of positive integers $\cW$, and 
sequences $\pmb{\varepsilon}=(\varepsilon_m)_{m\in \cW}$,
$\pmb{\delta}=(\delta_m)_{m\in \cW}$ of real numbers
satisfying
$0 \le \varepsilon_m < \delta_m \le 1$ for all $m \in \cW$,
we define
\begin{equation}\label{eqn:fJW}
\fJ(\cW,\pmb{\varepsilon},\pmb{\delta})=\bigcap_{m \in \cW} \fJ(m,\varepsilon_m,\delta_m) \, .
\end{equation}
To make the notation less cumbersome, we usually regard the values
$\varepsilon_m$ and $\delta_m$ as implicit,
and write $\fJ(m)$ for $\fJ(m,\varepsilon,\delta)$,
and $\fJ(\cW)$ for $\fJ(\cW,\pmb{\varepsilon},\pmb{\delta})$.

\begin{prop}\label{prop:sevenK}
Let $0<K_1 < K_2$ be real numbers. Let $\cW$ be a non-empty finite
set of integers such that every element $m \in \cW$
satisfies conditions (i)--(iii).
Suppose moreover, that for each $m \in \cW$, there
are real numbers $\varepsilon_m$, $\delta_m$ 
satisfying conditions (iv)--(vi).
Let $M=\lcm(\cW)$. Let $\fS \subset {[0,1]}$ be a finite
set of rational numbers $a/q$ (here $\gcd(a,q)=1$) 
with denominators $q$ bounded by $\sqrt[3]{M/2K_2}$.
Suppose that
\begin{equation}\label{eqn:decomp}
\fJ(\cW) \cap [0,M)
\; \subseteq \; \bigcup_{a/q \in \fS} 
\left( 
\frac{a}{q} M- \frac{\sqrt[3]{M/16}}{q} \; , \;
\frac{a}{q} M+ \frac{\sqrt[3]{M/16}}{q} 
\right).
\end{equation}
Then every odd integer $K_1 \le  N \le K_2$ is the sum of seven non-negative cubes.
\end{prop}
\begin{proof}
Let $N$ be an odd integer satisfying $K_1 \le N \le K_2$. 
It follows from assumptions (i)--(iii) that $M=\lcm(\cW)$ is squarefree
and divisible only by $3$ and primes $\equiv 5 \pmod{6}$. Thus there
exists a unique integer $T \in [0,M)$ such that
\begin{equation}\label{eqn:cuberoot}
N \equiv 8 T^3 \pmod{M}.
\end{equation}
Suppose $N$ is not the sum of seven cubes. Then, by Lemma~\ref{lem:sevenK},
for each $m \in \cW$, we have $\rem(T,m) \in [0,m) \setminus 
[\varepsilon_m \cdot m,\delta_m \cdot m)$. Thus $T \in \fJ(\cW) \cap [0,M)$.
By \eqref{eqn:decomp} there is some rational $a/q \in \fS$
such that
\[
- \frac{\sqrt[3]{M/16}}{q}< T-\frac{a}{q} M < \frac{\sqrt[3]{M/16}}{q},
\]
or equivalently
\[
-\frac{M}{2} < 8 (qT-aM)^3 < \frac{M}{2} \, .
\]
Moreover, the denominator $q$ is bounded by $\sqrt[3]{M/2K_2}$ and so
\[
q^3 N \le \frac{M N}{2K_2} \le \frac{M}{2}
\]
as $N \le K_2$. Hence 
\[
\lvert q^3 N - 8(qT-aM)^3 \rvert < M.
\]
However, by \eqref{eqn:cuberoot} we have
$q^3 N-8(qT-aM)^3 \equiv 0 \pmod{M}$.
Thus $q^3 N=8(aT-aM)^3$. It follows that $N$ is a perfect cube, and
so is certainly 
the sum of seven non-negative cubes.
\end{proof}

We shall mostly apply Proposition~\ref{prop:sevenK}
with the parameter choices given by the following lemma.
\begin{lem}\label{lem:parachoice}
Let $K\ge 10^5$. Let $K_1=9K/10$, $K_2=K$.
Let 
\begin{equation}\label{eqn:parainterval}
\frac{263}{100} K^{1/3} \le m \le \frac{292}{100} K^{1/3}. 
\end{equation}
Then
conditions (iv)--(vi) 
are satisfied with $\varepsilon_m=0$ and $\delta_m=1/10$.
\end{lem}

\section{Plan for the paper}\label{sec:plan}
The rest of the paper is devoted to understanding and computing
the intersections $\fJ(\cW) \cap [0,M)$ appearing in 
Proposition~\ref{prop:sevenK}. Section~\ref{sec:prop}
collects various properties of remainders and bad sets
that are used throughout. Section~\ref{sec:ripples}
provides justification, under a plausible assumption,
that the intersection $\fJ(\cW) \cap [0,M)$ should 
be decomposable as in \eqref{eqn:decomp}.
Section~\ref{sec:first} gives an algorithm (Algorithm~\ref{alg1})
which takes as
input a finite set of positive integers $\cW$ and an
interval $[A,B)$ and returns the intersection $\fJ(\cW) \cap [A,B)$.
We also give a heuristic analysis of the algorithm and its
running time. Section~\ref{sec:tower} introduces the 
concept of a \lq tower\rq, which a sequence 
\begin{equation}\label{eqn:tower}
\cW_0 \subseteq \cW_1 \subseteq \cW_2 \subseteq \cdots \subseteq \cW_r=\cW.
\end{equation}
Letting $M_i=\lcm(\cW_i)$, we prove the recursive
formula for computing $\fJ(\cW_i) \cap [0,M_i)$
in terms of $\fJ(\cW_{i-1}) \cap [0,M_{i-1})$.
This recursive formula together with Algorithm~\ref{alg1}
is the basis for a much more efficient algorithm (Algorithm~\ref{alg2})
for computing $\fJ(\cW) \cap [0,M)$ given in Section~\ref{sec:tower}.

In Section~\ref{sec:towercomp} we let
 $M^*$ be the product of all primes $p \le 167$
that are $\equiv 5 \pmod{6}$, and 
\begin{equation}\label{eqn:Wstar}
\cW^*=\{ \; m \mid M^* \quad : \quad 265 \times 10^{9} \le m \le
290 \times 10^{9}
\;
 \} \, .
\end{equation}
We use a tower and Algorithm~\ref{alg2} to compute  
$\fJ(\cW^*) \cap [0,M^*)$. The actual computation consumed 
about 18,300
hours of CPU time.

Section~\ref{sec:prooflarge} is devoted to proving Theorem~\ref{thm:main}
for $N\ge (9/10)^{3998}\cdot \cK
\approx 4.76 \times 10^{44}$, where $\cK=\exp(524)$. The approach is to
divide
the interval 
$(9/10)^{3998}  \cK \le N \le \cK$ into subintervals
$(9/10)^{n+1} \cK \le N \le (9/10)^n \cK$ with
$0 \le n \le 3997$, and apply Proposition~\ref{prop:sevenK}
and Lemma~\ref{lem:parachoice} to prove that all
odd integers in the interval $(9/10)^{n+1} \cK \le N \le (9/10)^n \cK$
are sums of seven non-negative cubes. 
Indeed, we show that given $0 \le n \le 3997$, there is some
suitable positive $\kappa$ such that the elements of $\cW_0=\kappa \cdot \cW^*$
satisfy conditions (i)--(iii) (with $K_1=(9/10)^{n+1} \cK$ 
and $K_2=(9/10)^n \cK$) and that moreover,
$\fJ(\cW_0)=\kappa \fJ(\cW^*)$. Thus the results of the huge computation of
Section~\ref{sec:towercomp} are recycled $3998$ times; on top of 
this $\cW_0$ we construct a tower and continue until we have found
a set $\cW$ that satisfies the hypotheses of Proposition~\ref{prop:sevenK},
thereby proving Theorem~\ref{thm:main} for $N\ge (9/10)^{3998}\cdot \cK$.
The CPU time for the computations described in Section~\ref{sec:prooflarge}
was around 10,000 hours.

The proof of Theorem~\ref{thm:main} is completed in 
Section~\ref{sec:proofsmall} where a  modified strategy is needed
to handle the \lq small\rq\ ranges $(9/10)^{n+1} \cK \le N \le (9/10)^n \cK$
with $3998 \le n \le 4226$. Although these intervals are small (and
few) compared to those handled in Section~\ref{sec:prooflarge},
we are unable to recycle
the computation of Section~\ref{sec:towercomp}. This makes the
computations far less efficient, though still practical.
The CPU time for the computations described in Section~\ref{sec:proofsmall}
was around 2,750 hours.

\section{Some Properties of Remainders and Bad Sets}\label{sec:prop}
\begin{lem}\label{lem:kappa}
Let $m$ and $\kappa$ be positive integers with $\kappa \mid m$. Then
for any real $x$ we have
\[
\quo\left(\frac{x}{\kappa},\frac{m}{\kappa} \right)=\quo(x,m),
\qquad
\rem\left(\frac{x}{\kappa},\frac{m}{\kappa} \right) =\frac{1}{\kappa} \rem(x,m).
\]
\end{lem}

Let $\kappa$ be a positive integer. For a set $X \subset \R$
we denote $\kappa X=\{ \kappa x : x \in X\}$. 

\begin{lem}\label{lem:kappa2}
Let $m$ and $\kappa$ be positive integers. 
Let $0 \le \varepsilon < \delta \le 1$ be real numbers.
Then
\[
\fJ(\kappa m, \varepsilon, \delta)=
\kappa \cdot \fJ(m,\varepsilon,\delta)\, .
\]
Let $\cW$ be a set of positive integers 
and for $m \in \cW$ let $0 \le \varepsilon_m < \delta_m \le 1$
be real numbers. Let 
\[
\cW^\prime=\kappa \cdot \cW \, ,
\quad
 \pmb{\varepsilon}=(\varepsilon_m)_{m \in \cW} \, , \quad
\pmb{\delta}=(\delta_m)_{m\in \cW} \, , \quad
\pmb{\varepsilon}^\prime=(\varepsilon_{m/\kappa})_{m \in \cW^\prime} \, \quad 
\pmb{\delta}^\prime=(\delta_{m/\kappa})_{m \in \cW^\prime} \, .
\]
Then
\[
\fJ(\cW^\prime, \pmb{\varepsilon}^\prime, \pmb{\delta}^\prime)
= \kappa \cdot \fJ(\cW,\pmb{\varepsilon},\pmb{\delta} ) \, .
\]
\end{lem}
\begin{proof}
By \eqref{eqn:fJ} and Lemma~\ref{lem:kappa},
\begin{alignat*}{2}
x \in \fJ(\kappa m, \varepsilon,\delta)
& \iff \rem(x,\kappa m) \in {[0, \kappa m)} \setminus {[\varepsilon \cdot \kappa m, \delta\cdot \kappa m)} \\
& \iff \frac{1}{\kappa} \rem(x,\kappa m) \in {[0,m)} \setminus {[\varepsilon\cdot m ,\delta\cdot m)}\\
& \iff \rem(x/\kappa,m) \in {[0,m)} \setminus {[\varepsilon \cdot m, \delta \cdot m)} \\
& \iff x/\kappa \in \fJ(m, \varepsilon, \delta) \\
& \iff x \in \kappa \cdot \fJ(m, \varepsilon,\delta). 
\end{alignat*}
This proves the first part of the lemma. The second part
now follows from   
\eqref{eqn:fJW}.
\end{proof}

\begin{lem}\label{lem:natural1}
Given positive integers $M_1 \mid M_2$, we define
the \lq natural\rq\ map
\[
\pi_{M_2,M_1} \; : \; {[0,M_2)} \longrightarrow {[0,M_1)},
\qquad x \mapsto \rem(x,M_1).
\]
Then
$\pi_{M_2,M_1}$ is surjective, and for any
$T \subseteq [0,M_1)$,
\[
\pi_{M_2,M_1}^{-1} (T)=\bigcup_{k=0}^{(M_2/M_1)-1} 
(k \cdot M_1 + T).
\]
\end{lem}
\begin{lem}
\label{lem:natural2}
Let $\cW_1$, $\cW_2$ be sets of positive integers with
$\cW_1 \subseteq \cW_2$. Write $\cU=\cW_2-\cW_1$. 
Let $M_i=\lcm(\cW_i)$. 
Write $\pi=\pi_{M_2,M_1}$.
Then $\pi(\fJ(\cW_2)) \subseteq \fJ(\cW_1)$ and
\[
\fJ(\cW_2) \cap [0,M_2)
={\pi^{-1} (\fJ(\cW_1) \cap [0,M_1))} \, \cap \,
\fJ(\cU) \, .
\]
\end{lem}
\begin{proof}
Let $x \in \R$ and let $y=\pi(x)=\rem(x,M_1)$. If $m \in \cW_1$
then $\rem(y,m)=\rem(x,m)$, as $m \mid M_1$. Observe that
\begin{alignat*}{2}
x \in \fJ(\cW_2)
& \iff 
\text{$\rem(x,m) \notin {[\varepsilon_m \cdot  m, \delta_m \cdot m)}$
for all $m \in \cW_2$} \\
& \implies 
\text{$\rem(x,m) \notin {[\varepsilon_m \cdot m, \delta_m \cdot m)}$
for all $m \in \cW_1$} \\
& \iff 
\text{$\rem(y,m) \notin {[\varepsilon_m \cdot m, \delta_m \cdot m)}$
for all $m \in \cW_1$} \\
& \iff y \in \fJ(\cW_1) \, .
\end{alignat*}
This shows that $\pi (\fJ(\cW_2)) \subseteq \fJ(\cW_1)$.
The rest of the lemma easily follows.
\end{proof}

\section{Gaps and Ripples}\label{sec:ripples}
We will soon give an algorithm for computing the intersection
\[
\fJ(\cW) \cap [0,M)=(\cap_{m \in \cW} \fJ(m)) \cap [0,M)
\qquad (M=\lcm(\cW))
\]
given a set $\cW$ that satisfies the conditions of 
Proposition~\ref{prop:sevenK}. The statement of Proposition~\ref{prop:sevenK}
(notably \eqref{eqn:decomp}) suggests that we are expecting
this intersection to be concentrated in small intervals around 
$aM/q$ for certain $a/q$ with relatively small denominators $q$.
In this section we provide an explanation for this. The situation
is easier to analyze if we make choices of parameters
as in Lemma~\ref{lem:parachoice}. Thus for this section we
fix the choices $\varepsilon_m=0$, $\delta_m=1/10$, and hence
$\fJ(m)=\fJ(m,0,1/10)$.
We suppose that the elements
$m \in \cW$ belong to an interval of the form
\begin{equation}\label{eqn:Linterval}
\frac{263}{100} L \le m \le \frac{292}{100} L,
\end{equation}
for some $L>0$ (\emph{c.f.} Lemma~\ref{lem:parachoice}). 
In fact, we show that if $q$ is large, and if the residues
of the integers $aM/m$ are regularly distributed modulo $q$
(in a sense that will be made precise), then the intersection
$\fJ(\cW) \cap [0,M)$ contains no points in a
certain explicitly given neighbourhood of $aM/q$. Likewise we show
for certain $a/q$ with $q$ small, that 
$\fJ(\cW) \cap [0,M)$ does contain some points
near $aM/q$. We stress that the material in this section
does not form part of our proof of Theorem~\ref{thm:main}.
It does however explain the results of our
computations that do form part of the proof of Theorem~\ref{thm:main},
and it lends credibility to them.

\medskip

We fix the following notation throughout this section.
\begin{itemize}
\item $L$ is a positive real number;
\item $\cW$ is a non-empty set of positive integers
that belong to the interval \eqref{eqn:Linterval};
\item $M=\lcm(\cW)$.
\end{itemize}

\subsection{Ripples}
\begin{prop}\label{prop:ripple1}
Suppose $M \ge 2000 L$.
Let $a/q \in [0,1)$ be a fraction in simplest form with
$1 \le q \le 9$ and $0 \le a \le q-1$. 
For $0 \le k \le 9-q$ let
\begin{equation}\label{eqn:psidef}
\psi_k = \frac{292}{100} \left(\frac{k}{q}+\frac{1}{10}\right),
\qquad
\Psi_k= \frac{263}{100} \cdot \frac{(k+1)}{q}.
\end{equation}
Then  $\psi_k < \Psi_k$ and
\begin{equation}\label{eqn:posrip}
\bigcup_{k=0}^{9-q} 
\left( 
\frac{a}{q} M + \psi_k \cdot L,
\frac{a}{q} M + \Psi_k \cdot L
\right) 
\subseteq \fJ(\cW) \cap [0,M).
\end{equation}
This recipe gives $103$ disjoint intervals contained in 
$\fJ(\cW) \cap [0,M)$
of total length $\xi \cdot L$ where
\[
\xi=\frac{261707}{10500} \approx 24.9 \qquad \text{(1 d.p.)}.
\]
\end{prop}
We shall informally refer to the union of intervals \eqref{eqn:posrip}
as \textbf{ripple emanating from $aM/q$} in the positive direction.
The reader will easily modify the proof below to show,
under similar hypotheses,
that there are ripples emanating from the $aM/q$ in the 
negative direction.
\begin{proof}
It is easy to check that $\psi_k<\Psi_k$ for $q \le 9$
and $0 \le k \le 9-q$.
The assumption $M \ge 2000 L$ ensures that the $103$ intervals
are contained in $[0,M)$ and are disjoint, so it is enough
to show that the intervals are contained 
in $\fJ(\cW)$.
Let $\alpha$ be a real number
belonging to the interval $\psi_k \cdot L < \alpha < \Psi_k \cdot L$.
We would like to show that $aM/q+\alpha \in \fJ(m)$ for
all $m \in \cW$. 
Let $m \in \cW$.
It follows from \eqref{eqn:Linterval} and \eqref{eqn:psidef} that
\begin{equation}\label{eqn:alpha}
\left(\frac{k}{q}  + \frac{1}{10}\right) m \le
\psi_k \cdot L  < \alpha < \Psi_k \cdot L  \le \frac{(k+1)}{q} m.
\end{equation}
As $m \mid M$ we can write
$aM=um$ with $u \in \Z$. Now $u=bq+s$ where $0 \le s \le q-1$.
Thus 
\[
\frac{a}{q}M=bm+\frac{s}{q} m.
\]
From \eqref{eqn:alpha},
\[
bm+\frac{(k+s)}{q} m+\frac{m}{10} < 
\frac{a}{q} M+\alpha < bm + \frac{(k+s+1)}{q} m \, .
\]
Let $k+s=qt+v$ where $0 \le v \le q-1$. 
Hence
\[
(b+t)m+\left(\frac{v}{q} +\frac{1}{10}\right) m < 
\frac{a}{q} M+\alpha 
< (b+t)m + \frac{(v+1)}{q} m \, .
\]
Observe that 
\[
\frac{1}{10} \le \frac{v}{q}+\frac{1}{10} < \frac{v+1}{q} \le 1,
\]
as $q \le 9$ and $0 \le v \le q-1$.
Thus $\quo(aM/q+\alpha,m)=b+t$ and
\[
\frac{m}{10}  
\; < \; \rem\left(\frac{aM}{q}+\alpha, m\right) 
\; < \; m.
\]
This shows that $aM/q+\alpha \in \fJ(m)$ as required.
\end{proof}

In the above proposition we showed the existence of ripples
emanating from $aM/q$ for $q \le 9$. There 
can also be ripples emanating for $aM/q$ for larger
values of $q$ if the sequence of residues $\overline{aM/m}$
in $\Z/q\Z$ contains large gaps as illustrated by the following
proposition.

\begin{prop}\label{prop:ripple2}
Let $a/q \in (0,1)$ be a rational in simplest form
with $q\ge 11$ and $1 \le a \le q-1$. 
Let $(q-10)/10<d < q-1$ be an integer, and let $s$ be a non-negative 
integer satisfying 
\begin{equation}\label{eqn:rineq}
s < q-d-1, \qquad
s < \frac{263}{290} (10d+10 -q).
\end{equation}
Suppose
\begin{equation}\label{eqn:aMm}
\overline{s+1}, \; \overline{s+2}, \dotsc, \overline{s+d}
\notin 
\left\{
\overline{aM/m} \; : \; m \in \cW
\right\}
\subseteq \Z/q\Z.
\end{equation}
Let
\[
\pi=\frac{292}{100}\cdot \frac{s}{q}, \qquad
\Pi= \frac{263}{100} \left( \frac{(s+d+1)}{q}  - \frac{1}{10} \right).
\]
Then $\pi < \Pi$ and
\[
\left( \frac{a}{q}M -\Pi \cdot L\,
, \, \frac{a}{q}M -\pi \cdot L \right)
\; \subseteq  \;
\fJ(\cW). 
\]
\end{prop}
\begin{proof}
Let $m \in \cW$, and recall that $m \mid M$.
Thus $aM/m$ is an integer, and hence so is $\rem(aM/m,q)$.
By assumption \eqref{eqn:aMm},
\[
\rem(aM/m,q) \ne s+1, s+2, \dotsc, s+d.
\]
Thus
$\rem(aM/m,q) \notin \left(s \, ,\, s+d+1 \right)$.
By Lemma~\ref{lem:kappa}, 
\[
\rem(aM/q,m) \, =\, \rem(aM,qm)\cdot 1/q \, =\, \rem(aM/m,q) \cdot m/q.
\]
Thus
\begin{equation}\label{eqn:remaM/q}
\rem(aM/q,m) \notin \left( \frac{s}{q} m \, , \, \frac{(s+d+1)}{q}  m \right).
\end{equation}
The condition $d> (q-10)/10$ implies that
\[
\frac{s}{q}  < \frac{(s+d+1)}{q}   - \frac{1}{10} \,.
\]
Let $\alpha$ belong to the interval 
\begin{equation}\label{eqn:alphainterval}
\frac{s}{q} m \; < \; \alpha \;  < \; \left( \frac{(s+d+1)}{q}   - \frac{1}{10} \right) m \,.
\end{equation}
We claim that
\[
\rem(aM/q-\alpha,m) \notin {[0,m/10)} .
\]
Suppose otherwise:
then we can write
\[
\frac{a}{q} M -\alpha=bm+r
\]
where $0 \le r < m/10$.
Thus
\[
bm+\frac{s}{q} m \, < \, bm+\alpha \, \le\, 
\frac{a}{q} M \, <\,  bm+\alpha+\frac{m}{10}
\, < \,
bm+\frac{(s+d+1)}{q} m
\]
as $\alpha$ satisfies \eqref{eqn:alphainterval}. This contradicts
\eqref{eqn:remaM/q}, and establishes our claim.
In fact we have shown that if $\alpha$ belongs to the interval
\eqref{eqn:alphainterval}, then $aM/q-\alpha \in \fJ(m)$.

\medskip

Suppose now that $\alpha$ belongs to the interval
$\pi \cdot L < \alpha < \Pi \cdot L$ (the second inequality 
in \eqref{eqn:rineq} ensures $\pi < \Pi$). To prove
the proposition, all we have to show is that $\alpha$
satisfies  the inequalities in \eqref{eqn:alphainterval}
for all $m \in \cW$. However, these follow straightforwardly
from the fact that 
all $m \in \cW$ belong to the interval~\eqref{eqn:Linterval}.
\end{proof}

A few remarks are in order concerning Proposition~\ref{prop:ripple2}
and its proof.
\begin{itemize}[leftmargin=5mm]
\item For simplicity we have only constructed the first interval
in a ripple emanating from $aM/q$ in the negative direction.
If inequalities \eqref{eqn:rineq} are satisfied with a
significant margin, then it is possible to construct
more intervals belonging to this ripple. 
Likewise, with a suitable modification of the assumptions
one can also construct a ripple in the positive direction.
\item The first inequality in \eqref{eqn:rineq} 
is imposed merely for simplicity;
if it does not hold one can also construct ripples emanating from $aM/q$
after suitably modifying the second inequality in \eqref{eqn:rineq}.
\item The one indispensable assumption in Proposition~\ref{prop:ripple2}
is the existence of a sequence 
\[
\overline{s+1}, \overline{s+2}, \dotsc, \overline{s+d}
\]
of consecutive residues belonging to 
$(\Z/q\Z)\setminus \{\overline{aM/m} : m \in \cW\}$
of length $d$ that is roughly larger than $q/10$. 
We shall show below that if there
is no such sequence, then $\fJ(\cW)$
contains no elements in a neighbourhood of $aM/q$.
\end{itemize}

\subsection{Gaps}

Let $a/q \in {[0,1]}$ be a rational in simplest form,
and let
\[
\Phi_{a/q} \; : \; \cW \rightarrow \Z/q\Z,
\qquad
m \mapsto \overline{a (M/m)}.
\]
In view of the above,
define the \textbf{defect} $d(\cW,a/q)$ of $\cW$ with respect to $a/q$
as the length of the longest sequence 
$\overline{s+1},\overline{s+2},\dotsc,
\overline{s+d}$ belonging to $(\Z/q\Z) \setminus \Phi_{a/q}(\cW)$.
As $\cW \ne \emptyset$, we have $d(\cW,a/q)<q$. For example, if 
$\Phi_{a/q}$ is surjective then $d(\cW,a/q)=0$, and if
$\Phi_{a/q}(\cW)=(\Z/q\Z)^*$ 
then $d(\cW,a/q)=1$. 

\begin{lem}\label{lem:ineq}
With notation as above, let $d=d(\cW,a/q)$. Let $x \in \R$.
Then there is some element $m \in \cW$ and an integer $k$ 
such that
\[
\left\lvert x-\frac{aM}{qm}-k \right\rvert \; \le \; \frac{d+1}{2q}.
\]
\end{lem}
\begin{proof}
Let $u \in \Z$ satisfy $\lvert u- qx \rvert \le 1/2$.
We first suppose that $d$ is even.
Consider the sequence
\[
\overline{u-d/2},\; \overline{u-d/2+1},\; \overline{u-d/2+2}, \dots,\; \overline{u+d/2}
\]
of $d+1$ elements of $\Z/q\Z$. By definition of $d$, one of these
equals $\Phi_{a/q}(m)$ for some $m \in \cW$.
Thus there is some integer $k$
such that
\[
\left\lvert u-\frac{aM}{m}-kq \right\rvert \; \le \;  \frac{d}{2}.
\]
As $\lvert u-qx \rvert \le 1/2$, the result follows.

Now suppose that $d$ is odd and $qx\ge u$ (the case $qx<u$
is similar). Consider the sequence
\[
\overline{u-(d-1)/2},\; \overline{u-(d-1)/2+1},\; \overline{u-(d-1)/2+2}, \dots,\; \overline{u+(d+1)/2}
\]
which again has $d+1$ elements, and so there is some $m \in \cW$
and some integer $k$ such that
\[
u-\frac{(d-1)}{2} \; \le \;  \frac{aM}{m}+kq \; \le \;  u+\frac{(d+1)}{2}.
\]
Since $0 \le qx -u \le 1/2$, the lemma follows.
\end{proof}

\begin{lem}\label{lem:m*}
Let 
\[
m^*=\frac{38398}{13875} \cdot L,
\]
Then for all $m \in \cW$,
\[
\left\lvert \frac{L}{m}-\frac{L}{m^*} \right\rvert \le
\frac{725}{38398}.
\]
\end{lem}
\begin{proof}
By \eqref{eqn:Linterval},
the quantity $L/m$ belongs to the interval
$[100/292, 100/263]$. We have chosen 
$m^*$ so that $L/m^*$
is the mid-point of the interval.
The lemma follows as $725/38398$ is half 
the length of the interval.
\end{proof}

\begin{prop}\label{prop:gaps}
With notation as above, let $d=d(\cW,a/q)$ and suppose that
$d< (q-10)/10$. Let 
\begin{equation}\label{eqn:mu}
\mu=\frac{38398}{725} \left(\frac{1}{20}-\frac{(d+1)}{2q}\right).
\end{equation}
Then
\[
\left(\frac{a}{q} M -\mu L\; ,\; \frac{a}{q}M+\mu L\right)
\; \cap \; \fJ(\cW)
=\emptyset.
\]
\end{prop}

A few words are perhaps appropriate to help the reader appreciate
the content of the proposition. We shall suppose that $q>11$.
If $\#\cW$ is large compared to $q$, then we expect that  
$\Phi_{a/q}$ is close to being surjective and which forces
 $d$ to be small.
If that is the case then $\mu$ should be close to 
$38398/(725 \times 20) \approx 2.64$. Suppose now that $\# \cW$
is large, but that $q$ is much larger.  
Suppose also that the residues in the image $\Phi_{a/q}(\cW)$
are \lq randomly\rq\ distributed in $\Z/q\Z$. The quantity $d$
measures how large the gaps between these residues in the image
can be, and we expect that $d$ should be around $q/\#\cW$.
We therefore
expect that $\mu \approx (38398/725)(\frac{1}{20}-\frac{1}{2 \cdot \#\cW})$.
We see that $\mu$ should be positive if $\cW$ has much more
than $10$ elements.

\begin{proof}[Proof of Proposition~\ref{prop:gaps}]
The assumption $d<(q-10)/10$ ensures that $\mu$
is positive. Let $y \in (aM/q-\mu L, aM/q+\mu L)$.
We would to like to show that there is some $m \in \cW$ such 
that $y \notin \fJ(m)$. 

Write $y=aM/q+\beta$ where
$\lvert \beta \rvert < \mu L$. Letting $x=1/20-\beta/m^*$
in Lemma~\ref{lem:ineq}, we deduce the existence of
some integer $k$ and some element $m \in \cW$ such that
\[
\left\lvert \frac{\beta}{m^*} + \frac{aM}{qm}+k -\frac{1}{20} \right\rvert \; \le \; \frac{(d+1)}{2q} \, .
\]
Thus
\[
\left\lvert \frac{\beta}{m} + \frac{aM}{qm}+k -\frac{1}{20} \right\rvert  \; \le \; 
\frac{(d+1)}{2q} + \left\lvert \frac{\beta}{m^*}-\frac{\beta}{m} \right\rvert \, .
\]
Using $\lvert \beta \rvert < \mu L$, Lemma~\ref{lem:m*}
and the definition of $\mu$ in \eqref{eqn:mu},
we see that 
\[
\left\lvert \frac{\beta}{m} + \frac{aM}{qm}+k -\frac{1}{20} \right\rvert  \; < \;
\frac{1}{20} \, .
\]
Thus $y=aM/q+\beta$ belongs to the interval 
$-km+(0,m/10)$, showing that $y \notin \fJ(m)$ as required.
\end{proof}

\section{A First Approach to Computing $\fJ(\cW)$}\label{sec:first}

In this section 
$\cW$ is a finite set of positive integers $m$.
Associated to each $m \in \cW$ are real numbers 
$0 \le \varepsilon_m < \delta_m<1$. We shall write
$\pmb{\varepsilon}=(\varepsilon_m)_{m\in \cU}$
and $\pmb{\delta}=(\delta_m)_{m \in \cU}$.

\begin{lem}\label{lem:ABintersect}
Let $A < B $ be real numbers. For $m \in \cW$, let 
\[
q_m=\quo(A,m), \qquad r_m=\rem(A,m).
\]
\begin{enumerate}
\item[(a)] Suppose $r_n \in {[\varepsilon_n \cdot n , \delta_n\cdot n)}$
for some $n \in \cW$. Write $A^\prime=\min((q_n+\delta_n)\cdot n, B)$.
Then
\[
\fJ(\cW)
\cap 
[A,B) 
=
\fJ(\cW) 
\cap  
[A^\prime,B) .
\]
\item[(b)] Suppose $r_m \notin {[\varepsilon_m \cdot m , \delta_m \cdot m)}$
for all $m \in \cW$.
Define
\begin{equation}\label{eqn:Aprime}
A_m=\begin{cases}
(q_m+\varepsilon_m)\cdot m & \text{if $r_m < \varepsilon_m \cdot m$}\\
(q_m+1+\varepsilon_m) \cdot m & \text{if $r_m \ge \delta_m \cdot m$},
\end{cases} \qquad
A^\prime=\min{\left(B,\min(A_m)_{m \in \cW}\right)}.
\end{equation}
Then
\[
\fJ(\cW) \cap [A,B) \;
=
\;
\left( 
\fJ(\cW) 
\cap  
[A^\prime,B) 
\right)
\; \cup \; 
[A,A^\prime) 
.
\]
\end{enumerate}
\end{lem}
\begin{proof}
Suppose $n \in \cW$ satisfies
$r_n \in [\varepsilon_n \cdot n, \delta_n \cdot n)$, and let $A^\prime$
be as in (a).
By \eqref{eqn:fJ} we have
\[
(q_n \cdot n+[\varepsilon_n \cdot n, \delta_n \cdot n)) \cap \fJ(n)
=\emptyset.
\]
Observe that 
$[A,A^\prime) \subseteq q_n \cdot n+[\varepsilon_n \cdot n, \delta_n \cdot n)$ 
and $[A,A^\prime) \subseteq [A,B)$. Part (a) follows.

\medskip

Suppose now that $r_m \notin [\varepsilon_m \cdot m, \delta_m \cdot m)$
for all $m \in \cW$, and let $A^\prime$ be as in (b). 
It is easy to check that $\rem(A^{\prime\prime},m) \notin [\varepsilon_m \cdot m,
\delta_m \cdot m)$ for all $A^{\prime\prime} \in [A,A_m)$.
From this we see that 
$[A,A^\prime) \subseteq \bigcap_{m \in \cW} \fJ(m,\varepsilon_m,\delta_m)
=\fJ(\cW)$.
Part (b) follows.
\end{proof}

Lemma~\ref{lem:ABintersect} immediately leads us to the following 
algorithm.

\begin{alg}\label{alg1}
To compute 
$ \fJ(\cW) \cap {[A,B)}$
as a disjoint union of intervals $\bigcup_{I \in \cI} I$.

\noindent \textbf{Input:} $A$, $B$, $\cW$, $\pmb{\varepsilon}$,
$\pmb{\delta}$.

\noindent Initialize $\cI \leftarrow \emptyset$.

\noindent Repeat the following steps until $A=B$.
\begin{enumerate}
\item[(a)] Loop through the elements $m \in \cW$ computing 
$q_m=\quo(A,m)$, $r_m=\rem(A,m)$.
\item[(b)] If there is some $n \in \cW$ such that 
$\varepsilon_n \cdot n \le r_n <
\delta_n \cdot n$ then 
\[
A \leftarrow \min\left((q_n +\delta_n) \cdot n,B\right)
\]
and go back to (a).
\item[(c)] Otherwise, let $A^\prime$ be as in \eqref{eqn:Aprime}.
Let $\cI\leftarrow \cI \cup \{[A,A^\prime)\}$ and then $A \leftarrow A^\prime$.
Go back to (a).
\end{enumerate}
\noindent \textbf{Output:} $\cI$.
\end{alg}

\subsection*{A heuristic analysis of Algorithm~\ref{alg1}
and its running time}

Let $x \in [0,M)$ and recall that $\rem(x,m) \in [0,m)$. Moreover, 
$x \in \fJ(m)$
if and only if 
$\rem(x,m) \in [0,m) \setminus [\varepsilon_m \cdot m,\delta_m \cdot m)$. 
Thus the \lq probability\rq\ that $x$ 
belongs to $\fJ(m)$ is $1-(\delta_m-\varepsilon_m)$. 
Assuming \lq independence of events\rq\ 
we expect that the 
total length of intervals produced by Algorithm~\ref{alg1}
is
\begin{equation}\label{eqn:outputlength}
(B-A) \cdot \prod_{m \in \cW} (1-\delta_m+\varepsilon_m) \, .
\end{equation}

To analyse the running time, we shall suppose parameter
choices as in Lemma~\ref{lem:parachoice}: namely
$\varepsilon_m=0$ and $\delta_m=1/10$ for all $m \in \cW$.
Moreover, we shall suppose that the elements of $m\in \cW$
belong to an interval \eqref{eqn:Linterval} for some large positive $L$.
By the above, the expected
total length of the intervals produced by Algorithm~\ref{alg1}
is $(B-A) \cdot 0.9^{\# \cW}$.
Moreover, we suppose that $\cW$ is sufficiently large so that
the length of the output should be negligible compared to $B-A$;
this should mean that step (c) is relatively rare.
We will estimate
the expected number of times we loop through steps (a), (b). 
Note that in step (b), $A$ is increased by $0.1\cdot n-r_n$.
The remainder $r_n=\rem(A,n)$ belongs to $[0,0.1 \cdot n)$. 
We regard the increase as a product $(0.1-r_n/n) \cdot n$.
Treating $r_n/n$ as a random variable uniformly distributed
in $[0,0.1)$ and $n$ as a random variable uniformly distributed
in interval \eqref{eqn:Linterval}, we see that the expected
increase is $0.05 \cdot (2.63+2.92)L/2=0.13875\cdot L$.
A standard probability theory argument that we omit tells us that the
expected number of times the algorithm loops through steps (a), (b)
is roughly
\[
(B-A)/(0.13875 L) \approx 7(B-A)/L \, .
\]

We now suppose that $K$ is very large, and we would like
to compute the intersection
$\fJ(\cW) \cap [0,M)$ for some set $\cW$ where
we hope that the hypotheses of 
Proposition~\ref{prop:sevenK} and
Lemma~\ref{lem:parachoice} 
are satisfied.  In particular,
we take $L=K^{1/3}$.
The number of steps should be around
$7 M/K^{1/3}$.
We have to choose $\cW$ so that $M=\lcm(\cW)$ is much
larger than $K$ (see \eqref{eqn:decomp} and just above it). 
Thus the number of steps
to compute $\fJ(\cW)$ is
much greater than $K^{2/3}$. 
For $K=\exp(524)$, the expected number of steps
is larger than $10^{150}$, which makes the computation
entirely impractical.

\section{A Refined Approach to Computing $\fJ(\cW)$: The Tower}
\label{sec:tower}

In this section we let $\cW$ be a set of positive integers
with $M=\lcm(\cW)$. Let $M_0, M_1, M_2, \hdots, M_r$
be positive integers such that $M_i \mid M_{i+1}$ and $M_r=M$.
Write $p_i=M_{i+1}/M_i$. In our later computations the
$p_i$ will be primes, but we need not assume that yet.
Let
\[
\cW_i=\{ m \in \cW \; : \; m \mid M_i \}.
\]
We suppose that $M_i=\lcm(\cW_i)$. Write $\cU_i=\cW_{i+1} \setminus \cW_i$.
Recall (Lemmas~\ref{lem:natural1} and~\ref{lem:natural2}) 
that we have natural surjections $\pi_{j,i} \; : \; [0,M_j) \rightarrow [0,M_i)$
whenever $j\ge i$, and that these restrict to give maps (not necessarily
surjections) 
$\fJ(\cW_j) \rightarrow \fJ(\cW_i)$.
We shall refer to the sequence of inclusions \eqref{eqn:tower}
as
a \textbf{tower leading up to} $\fJ(\cW)$, and use this
to compute $\fJ(\cW)$.

\begin{lem} \label{lem:pullback}
Let $0 \le i \le r-1$.
Suppose $\cI_i$ is a finite set of disjoint subintervals 
of $[0,M_i)$ such that
\[
\fJ(\cW_i) \cap [0,M_i) = 
\bigcup_{I \in \cI_i} I.
\]
Then
\[
\fJ(\cW_{i+1}) \cap [0,M_{i+1})
=\bigcup_{I \in \cI_i}
\bigcup_{k=0}^{p_i -1} 
\left((k \cdot M_i+I) \cap 
\fJ(\cU_i)
\right)
\, .
\]
\end{lem}
\begin{proof}
This is immediate from Lemmas~\ref{lem:natural1} and~\ref{lem:natural2}.
\end{proof}



Lemma~\ref{lem:pullback} immediately leads us to the following algorithm.
\begin{alg}\label{alg2}
The following computes 
a finite set $\cI=\cI_r$ of subintervals
of $[0,M)$ such that
$\fJ(\cW) \cap [0,M)= \bigcup_{I \in \cI} I$.

\noindent \textbf{Input:} $\cW_0, \dotsc, \cW_r=\cW$,
$\pmb{\varepsilon}$,
$\pmb{\delta}$.

\noindent Initialize: $\cI_0$ to be the set of 
disjoint intervals  whose union equals $\fJ(\cW_0) \cap [0,M_0)$,
and which is computed using Algorithm~\ref{alg1}.

\noindent Initialize: $i \leftarrow 0$.

\noindent Repeat the following steps until $i=r$.
\begin{enumerate}
\item[(a)] $\cI_{i+1} \leftarrow \emptyset$.
\item[(b)] for $I \in \cI_i$ and $k \in \{0,\dotsc,p_i-1\}$,
compute, using Algorithm~\ref{alg1}, a finite set $\cI^\prime$ 
of subintervals of $[0,M_{i+1})$ such that 
$(k \cdot M_i + I) \cap \fJ(\cU_i)= \bigcup_{I^\prime \in \cI^\prime} I^\prime$; let $\cI_{i+1} \leftarrow \cI_{i+1} \cup \cI^\prime$. 
\item[(c)] $i \leftarrow i+1$.
\end{enumerate}
\noindent \textbf{Output:} $\cI=\cI_r$.
\end{alg}

\subsection*{A heuristic analysis of Algorithm~\ref{alg2}
and its running time}

We shall suppose, as in Lemma~\ref{lem:parachoice},
 that 
$\varepsilon_m=0$ and $\delta_m=1/10$ for all $m \in \cW$.
Write $n_i=\# \cW_i$. We 
assume that the elements of $\cW_i$, $\cU_i$ belong to an interval of the form 
$[263 L/100, 292 L/100]$ for some large $L$.
 By our previous analysis, we expect
that we can compute $\cI_0$ in roughly $7 M_0/ L$
steps.
The total length $\ell(\cI_0)$ of the intervals in $\cI_0$
should roughly be $0.9^{n_0} M_0$. In Step (b) of the algorithm,
we will replace each $I \in \cI_0$ with $p_0$ intervals of the same
length, and then apply Algorithm~\ref{alg1} to each. Thus
we expect that the number of steps to compute $\cI_1$ to be
roughly
\[
\frac{7p_0 \cdot 0.9^{n_0} \cdot M_0}{L}\approx 
\frac{7 M_1 \cdot 0.9^{n_0}}{L}.
\] 
The total length of the intervals in $\cI_1$ should 
be roughly $M_1 \cdot 0.9^{n_1}$. It is now apparent that
the total number of steps should be around
\[
(7/L) \cdot (M_0+M_1 \cdot 0.9^{n_0}
+M_2 \cdot 0.9^{n_1}+\cdots +M_r \cdot 0.9^{n_{r-1}}) \, .
\]



\section{A Large Computation}\label{sec:towercomp}
Let $M^*$ be the product of all primes $p \le 167$
that are $\equiv 5 \pmod{6}$, and $\cW^*$ is as in \eqref{eqn:Wstar}.
In this section we compute $\fJ(\cW^*) \cap [0,M^*)$, using
a tower and Algorithm~\ref{alg2}. 
As explained in the plan (Section~\ref{sec:plan}),
the result of this computation will be
reused again and again in Section~\ref{sec:prooflarge}. 
Let
\[
M_0=5\times 11 
\times 17
\times 23
\times 
29
\times
41
\times
47
\times 53
\times 
59
\times
71
\times
83
\times
89 \, ,
\]
which is the product of the primes $<100$ that are $\equiv 5 \pmod{6}$. Let
\begin{gather*}
M_1=101 \cdot M_0,
\quad
M_2=107 \cdot M_1,
\quad
M_3=113 \cdot M_2,
\quad
M_4=131 \cdot M_3,\\
M_5=137 \cdot M_4,
\quad
M_6=149 \cdot M_5,
\quad
M^*=M_7=167 \cdot M_6.
\end{gather*}
We let
\[
\cW_i=\{ \; m \mid M_i \quad : \quad 265 \times 10^{9} \le m \le
290 \times 10^{9}
\;
 \}. 
\]
Thus $\cW_0 \subseteq \cdots \subseteq \cW_7=\cW^*$. 
We checked that $M_i=\lcm(\cW_i)$. 
Table~\ref{table1} gives the cardinalities of the $\cW_i$. 
We use this tower and Algorithm~\ref{alg2} to compute $\fJ(\cW^*) \cap [0,M^*)$.
By our heuristic in the previous section, 
the number of steps needed for this computation
should very roughly be equal to $6.0 \times 10^{10}$, which is the sum
of the entries of the table's third column.
It appears from this estimate
that the computation can be done in reasonable time.

\begin{table}
\caption{The $M_i$ and the $\cW_i$ are given at the beginning of Section~\ref{sec:towercomp}.
The third column gives an estimate for the number of steps needed to compute
$\fJ(\cW_i) \cap [0,M_i)$ from $\fJ(\cW_{i-1}) \cap [0,M_{i-1})$
according to the heuristic analysis at the end of Section~\ref{sec:tower}.}
\label{table1}
\begin{centering}
{\tabulinesep=1.2mm
\begin{tabu}{|c|c|c|c|}
\hline\hline
\rule[-2.0ex]{0pt}{5ex}
$i$
& $\log_{10}(M_i)$ (1 d.p.) & 
$\displaystyle n_i=\# \cW_i$ & $\frac{7 M_i \cdot 0.9^{n_{i-1}}}{10^{11}}$ (2 s.f.)\\
\hline\hline
$0$ 
& $18.3$ & $16$ & $1.4 \times 10^9$\\
\hline
$1$ 
&  $20.3$ & $38$ & $2.6 \times 10^9 $ \\
\hline
$2$ 
& $22.3$ & $83$ & $2.7 \times 10^{10}$\\
\hline
$3$ 
& $24.4$ & $149$ & $ 2.7 \times 10^{10}$\\
\hline
$4$ 
& $26.5$ & $250$ & $3.4 \times 10^{9}$\\
\hline
$5$ 
& $28.6$ & $401$ & $1.1 \times 10^7 $\\
\hline
$6$ 
& $30.8$ & $620$ & $2.0 \times 10^2$\\
\hline
$7$ 
& $33.0$ & $911$ & $3.2 \times 10^{-6}$\\
\hline
\end{tabu}
}
\end{centering}
\end{table}

We wrote simple implementations of Algorithms~\ref{alg1}
and~\ref{alg2} for the computer algebra system \texttt{Magma} \citep{magma}.
We divided the interval $[0,M_0)$ into $59,000$ 
subintervals of equal length and ran our program
on each of these intervals $[A_{k-1},A_k)$ 
successively computing 
$ \fJ(\cW_i) \cap \pi_{i,0}^{-1} ([A_{k-1},A_{k}))$
for $i=0,\dotsc,7$.  
Our computation was distributed over 59 processors
(on a 64 core machine with 2500MHz AMD Opteron Processors).
Note that 
\[
\fJ(\cW_i) \cap [0,M_i) = \bigcup_{k=1}^{59,000} 
\fJ(\cW_i)
\cap 
\pi_{i,0}^{-1} ([A_{k-1},A_{k})) 
;
\]
thus our computation gives us a decomposition of $\fJ(\cW_i) \cap [0,M_i)$ as a union
of disjoint intervals.  The total CPU time for the computation was
around 18,300 hours, but as we distributed the computation over
 59 processors, it was over in less than two weeks. 

\begin{table}
\caption{Some details for the computation described Section~\ref{sec:towercomp}.
The second column gives $\# \cI_i$, where $\cI_i$ is a disjoint collection
of intervals such $\cup \cI_i=\fJ(\cW_i) \cap [0,M_i)$. The third column gives
the total length $\ell_i$ of these intervals. The fourth column gives
(2 s.f.) the ratio $\ell_i/M_i$. According to the heuristic
at the end of Section~\ref{sec:first}, this ratio should approximately equal $0.9^{n_i}$
which is given (2 s.f.) in the last column (here $n_i=\# \cW_i$
as in Table~\ref{table1}). We explain the discrepancy between the last two columns
in Subsection~\ref{sub:sanity}.}
\label{table2}
\begin{centering}
{\tabulinesep=1.2mm
\begin{tabu}{|c|c|c|c|c|}
\hline\hline
$i$
& $\# \cI_i$ & $\ell_i=\ell(\fJ(\cW_i) \cap [0,M_i))$ & 
$\ell_i/M_i$ &
$0.9^{n_i}$
\\
\hline\hline
$0$ &
$23,458,002$
& $365,300,497,739,376,385 \frac{8}{10}$
& $1.85 \times 10^{-1}$ & $1.85 \times 10^{-1}$ 
\\
\hline
$1$ &
$553,209,618$ & 
$3,625,384,986,862,035,664\frac{4}{10}$ &
$1.82 \times 10^{-2}$ & $1.82 \times 10^{-2}$ \\
\hline
$2$ &
$1,106,375,245$
& 
$3,313,998,145,602,553,709\frac{1}{10}$
& $1.56 \times 10^{-4}$
& $1.59 \times 10^{-4}$
\\ 
\hline
$3$
&
$209,982,392$
&
$350,826,426,611,537,217\frac{1}{10}$
&
$1.46\times 10^{-7}$
&
$ 1.52\times 10^{-7} $
\\ 
\hline
$4$
&
$1,062,201$
& 
$1,076,402,154,947,217\frac{8}{10}$ &
$3.41\times 10^{-12}$ & 
$3.64\times 10^{-12}$
\\ 
\hline
$5$
&
$904$ 
&
$20,663,973,893,432\frac{1}{10}$
&
$4.78\times 10^{-16}$ 
&
$4.48\times 10^{-19}$
\\ 
\hline
$6$
&
$870$
& 
$20,504,346,087,851\frac{7}{10}$ 
&
$3.19 \times 10^{-18}$ 
&
$4.27\times 10^{-29}$
\\ 
\hline
$7$
&
$861$
&
$20,382,195,221,000\frac{6}{10}$
& 
$1.90\times 10^{-20}$ 
& 
$2.07\times 10^{-42}$
\\ 
\hline
\end{tabu}
}
\end{centering}
\end{table}

\begin{lem}\label{lem:Mprime}
There are 
sequences $(B_j)_{j=1}^{854}$, $(C_j)_{j=1}^{854}$
 contained in $[0,M^*]$ such that
\[
B_1<C_1<B_2<C_2<\cdots<B_{854}<C_{854}
\]
and
\[
\fJ(\cW^*) \cap [0,M^*)=\bigcup_{j=1}^{854} [B_j,C_{j}),
\] 
with total length 
$\displaystyle \sum_{j=1}^{854} (C_{j}-B_j)=
20382195221000\frac{6}{10}$.
\end{lem}
\begin{proof}
As indicated by Table~\ref{table2},
our computation gives $\fJ(\cW^*) \cap [0,M^*)$
as a union of $861$ intervals disjoint subintervals of $[0,M^*)$.
Among these there are $7$ pairs of the form 
$[\alpha,\beta) \cup [\beta,\gamma)$, where the values of $\beta$
are of the form $\beta^\prime \cdot M^*/59000$ with
\[
\beta^\prime=7375, \quad
14750, \quad
22125, \quad
29500, \quad
36875, \quad
44250, \quad
51625 \, .
\]
These subdivisions are clearly a result of our original subdivision
of interval $[0,M_0^*)$ into $59000$ subintervals of equal length.
We simply replace the pairs $[\alpha,\beta) \cup [\beta,\gamma)$
with $[\alpha,\gamma)$ so that $\fJ(\cW^*) \cap [0,M^*)$ is expressed
as a union of $854$ intervals. This simplification of course
preserves the total length of intervals.
\end{proof}

\subsection{Remarks and Sanity Checks}\label{sub:sanity}
Our computations 
are done with exact arithmetic. The reader
will note by looking back at Algorithms~\ref{alg1} and~\ref{alg2}
(and recalling that all $\varepsilon_m=0$ and $\delta_m=m/10$)
that the end points of the intervals encountered will be
rationals with denominators that are divisors of $10$,
except for the $A_k$ appearing in our original subdivision which have
denominators that are divisors of $59,000$.
As a check on our computations, we verify that our results
for $\fJ(\cW^*) \cap [0,M^*)$ 
are consistent with Proposition~\ref{prop:ripple1}.
The set $\cW^*$ satisfies
\[
\min(\cW^*)=265,024,970,473 \qquad
\max(\cW^*)=289,916,573,827.
\]
We take $L=
\min(\cW^*) \cdot 100/263$.
It turns out that
$L >\max(\cW^*) \cdot  100/292$.
Thus $\cW^*$ is contained in the interval~\eqref{eqn:Linterval}
for this value of $L$.
Proposition~\ref{prop:ripple1} yields a total $103$ intervals of the 
form $(aM^*/q+\psi_k \cdot L,aM^*/q+\Psi_k\cdot L)$ that must be
contained in $\fJ(\cW^*) \cap [0,M^*)$. We checked that each of these
is contained in one of the $854$ intervals produced by our 
computation. 
It is instructive to compare the fourth and fifth columns of Table~\ref{table2}.
According to our heuristic, the total length $\ell(\fJ(\cW_i) \cap [0,M_i))$ 
should be around $M_i \cdot 0.9^{n_i}$ (with $n_i=\# \cW_i$) and therefore we expect
the two columns to be roughly the same. 
From the table, we see that this heuristic
is remarkably accurate for $0 \le i \le 4$, and extremely inaccurate for $i \ge 5$.
An explanation for this is provided by the ripples.
The total length of the intervals contained in
$\fJ(\cW_i) \cap [0,M_i)$ 
produced by
 Proposition~\ref{prop:ripple1} is $\approx 24.9 L$. 
Now 
\[
\frac{24.9 L}{M_5}=5.8 \times 10^{-17}, \qquad
\frac{24.9 L}{M_6}=4.0 \times 10^{-19}, \qquad
\frac{24.9 L}{M_7}=2.3 \times 10^{-21} \qquad \text{($1$ s.f.)},
\]
which does provide an explanation for the discrepancy between the too columns. 
Proposition~\ref{prop:ripple2} (with $\cW^*$ and $M^*$ in place of $\cW$ and $M$)
produces $172$ intervals with $11 \le q \le 100$ with total length $\approx 17.8 L$.
We checked that each of these is also contained in one of the $854$ intervals
produced by our computation.

\medskip

According to the overall philosophy of Section~\ref{sec:ripples}, 
the set
$\fJ(\cW^*) \cap [0,M^*)$ should be concentrated in short intervals 
around rational multiples $(a/q) \cdot M^*$
with $q$ small. To test this,
we computed, using continued fractions, the best
rational approximation to $(B_i+C_i)/(2M^*)$ 
with denominator at most $10^{20}$,
for $1 \le i \le 854$. The largest denominator we found was $42$.

\medskip

The reader is probably wondering,
given that we are employing $59$ processors, why we have subdivided
$[0,M_0)$ into $59,000$ intervals instead of $59$ intervals. 
This was done purely for memory management reasons. A glance
at Table~\ref{table2} will show the reader that there is an explosion
of intervals at levels $i=1$, $2$, $3$. By dividing $[0,M_0)$ into
$59,000$ subintervals, we only need to store roughly $1/59000$-th
of the intervals appearing at levels $i$ at any one time per processor,
and so only need to store around $1/1000$-th of these intervals
in the memory at any one time. 

\section{Proof of Theorem~\ref{thm:main} for $N\ge (9/10)^{3998}\cdot \exp(524)
\approx 4.76 \times 10^{44}$}\label{sec:prooflarge}

The reader might at this point find it helpful
to review the first paragraph 
of Section~\ref{sec:criterion} as well as the plan in Section~\ref{sec:plan}.
Let $\cK=\exp(524)$.
In this section we prove Theorem~\ref{thm:main} 
for $N \ge (9/10)^{3998} \cK$.
We shall divide the interval 
$(9/10)^{3998}  \cK \le N \le \cK$ into subintervals
$(9/10)^{n+1} \cK \le N \le (9/10)^n \cK$ with
$0 \le n \le 3997$. We apply Proposition~\ref{prop:sevenK}
and Lemma~\ref{lem:parachoice} to prove that all
odd integers in the interval $(9/10)^{n+1} \cK \le N \le (9/10)^n \cK$
are sums of seven non-negative cubes.

\begin{lem}\label{lem:kappadef}
Let $0 \le n \le 3997$. Let $K=(9/10)^n \cdot \cK$.
There exists an integer $\kappa$ that satisfies
\begin{enumerate}
\item[(a)] $\kappa$ is squarefree;
\item[(b)] $3 \mid \kappa$
\item[(c)] $\kappa/3$ is divisible only by primes $q \equiv 5 \pmod{6}$
that satisfy $q>167$;
\item[(d)] $\kappa$ belongs to the interval
\begin{equation}\label{eqn:kappainterval}
\frac{263}{265} \cdot \frac{K^{1/3}}{10^{11}} \;
\le \; \kappa \;
\le \;
\frac{292}{290} \cdot \frac{K^{1/3}}{10^{11}} \, . 
\end{equation}
\end{enumerate}
\end{lem}
\begin{proof}
We proved the lemma using a \texttt{Magma} script.
Let $I_1$, $I_2$ be the lower and upper bounds for $\kappa$
in \eqref{eqn:kappainterval}. If $I_2<10^7$ then our script
uses brute enumeration of integers in the interval $[I_1,I_2]$
to find a suitable $\kappa$. Otherwise, the
script takes $\tau$
to be a product of consecutive primes $\equiv 5 \pmod{6}$
starting with $173$ up to a certain bound, and keeps increasing
the bound until $I_2/\tau<10^7$. It then loops through
the integers $I_1/3 \tau \le \mu \le I_2/3 \tau$
until it finds one such that $\kappa=3 \mu \tau$
satisfies conditions (a), (b), (c).
\end{proof}
\noindent \textbf{Remark.}
For $n=3998$, the interval in \eqref{eqn:kappainterval} is 
$7481.6\dots \le \kappa \le 7590.5\dots$, which is too short for the 
existence a suitable $\kappa$. This is also the case for
most values
of $n$ that are $\ge 3998$.

\begin{lem}\label{lem:bootstrap}
Let $0 \le n \le 3997$ and let $\kappa$
be as in Lemma~\ref{lem:kappadef}. 
 Let $\cW^*$ and $M^*$
be as in Lemma~\ref{lem:Mprime}.
Let
\[
\cW_0=\{\kappa \cdot m^* \; : \; m^* \in \cW^* \}, \qquad M_0=\lcm(\cW_0)=\kappa M^* .
\]
Let $\varepsilon_m=0$ and $\delta_m=1/10$ for all $m \in \cW_0$.
Then $m \in \cW_0$ satisfy the conditions (i)--(vi)
of Section~\ref{sec:criterion}, where
\[
K_1=(9/10)^{n+1} \cdot \cK, \qquad K_2=(9/10)^n \cdot \cK \, .
\]
Moreover,
\begin{equation}\label{eqn:fJcW0}
\fJ(\cW_0) \cap [0,M_0) \, = \, 
\bigcup_{j=1}^{854} [ \kappa \cdot B_j \, , \, \kappa \cdot C_{j}),
\end{equation}
where the $B_j$ and $C_j$ are as in Lemma~\ref{lem:Mprime}. 
\end{lem}
\begin{proof}
All $m^* \in \cW^*$ are squarefree and divisible only 
by primes $q\le 167$ satisfying $q \equiv 5 \pmod{6}$.
Thus conditions (i)--(iii) of Section~\ref{sec:criterion} 
are satisfied by $m \in \cW_0$. As we are taking
$\varepsilon_m=0$ and $\delta_m=1/10$, to verify
conditions (iv)--(vi) 
we may apply Lemma~\ref{lem:parachoice}.
For this we need only check \eqref{eqn:parainterval}
holds for $m \in \cW_0$, where $K=K_2$. This immediately
follows from \eqref{eqn:kappainterval} and the fact
that $\cW^* \subset [265 \times 10^9, 290 \times 10^9]$.

Finally, by Lemma~\ref{lem:kappa2},
\[
\fJ(\cW_0) \cap [0,M_0) \, = \, 
\kappa \cdot (\fJ(\cW^*) \cap [0,M^*)).
\]
Lemma~\ref{lem:Mprime} completes the proof.
\end{proof}

Our \texttt{Magma} script for proving Theorem~\ref{thm:main} in the range
$K_1 \le N \le K_2$ proceeds as follows. We inductively construct
a tower $W_0 \subset W_1 \subset W_2 \subset \cdots$.
Observe that
\[
\frac{\ell(\fJ(\cW_0) \cap [0,M_0))}{M_0} =\frac{\ell(\fJ(\cW^*) \cap [0,M^*))}{M^*}
\approx 1.90 \times 10^{-20},
\]
thus the computation of the previous section has already substantially depleted
the interval $[0,M_0)$. 
Given $\cW_{i}$, and $M_{i}$, we let
$p_{i}$ be the smallest prime $\equiv 5 \pmod{6}$ that does
not divide $M_{i}$ and let $M_{i+1}=p_{i} M_{i}$.
The script then writes down positive integers $m$ belonging to
the interval \eqref{eqn:parainterval}, such 
that $m \mid M_{i+1}$ and $3 p_{i} \mid m$. It is not
necessary or practical to find all such integers, but
we content ourselves with finding around 
$3 \log(p_{i})/\log(0.9^{-1})$
of them; we explain this choice shortly. 
These $m$ will form the set $\cU_{i}$ and we take
$\cW_{i+1}=\cW_{i} \cup \cU_{i}$. The script then applies
our implementation of Algorithm~\ref{alg2} to
compute $\fJ(\cW_{i+1}) \cap [0,M_{i+1})$ as a union of disjoint intervals.
Our heuristic analysis of Algorithm~\ref{alg2}
suggests that $\ell(\fJ(\cW_{i+1}) \cap [0,M_{i+1}))$ should roughly equal
$p_{i} \cdot 0.9^{\#\cU_{i}} \cdot \ell(\fJ(\cW_{i}) \cap [0,M_i))$.
We desire the total length of the intervals to decrease in each
step of the tower, so
we should require $\# \cU_{i} > \log(p_{i})/\log(0.9^{-1})$. 
Experimentation suggests that requiring $\# \cU_{i} \approx 3 \log(p_{i})/\log(0.9^{-1})$ 
provides good control of both the total length of $\fJ(\cW_{i}) \cap [0,M_i)$
and the number of intervals that make it up. Our script continues to
build the tower and compute successive $\fJ(\cW_i) \cap [0,M_i)$ until it finds $\cW=\cW_i$
and $M=M_i$ that satisfy \eqref{eqn:decomp} for some set of rationals $\fS \subset [0,1]$
with denominators bounded by $\sqrt[3]{M/2K}$. Specifically, once $M_i>2K$, 
for each of the disjoint intervals $[\alpha,\beta)$ 
that make up $\fJ(\cW) \cap [0,M)$, 
the script
uses continued fractions to compute the best rational approximation $a/q$
to $(\alpha+\beta)/2M$ with $q \le \sqrt[3]{M/2K}$, and then checks whether
$[\alpha,\beta) \subseteq (aM/q - \sqrt[3]{M/16}/q \, , \, aM/q + \sqrt[3]{M/16}/q )$.
The script continues constructing the tower until this criterion is satisfied
for all the intervals making up $\fJ(\cW)$. It then follows from
Proposition~\ref{prop:sevenK} that all odd integers in the range
$\cK \cdot (9/10)^{n+1} \le N \le \cK \cdot (9/10)^n$ are sums of seven non-negative
cubes. We again distributed the computation among 59 processors on the afore-mentioned machine, 
with each processor
handling an appropriate portion of the range $0 \le n \le 3997$. The script succeeded
in finding an appropriate $\cW$ for all $n$ in this range. The entire CPU time
was around 10,000 hours, but as the computation was distributed among 
59 processors the actual time was around 7 days.

\bigskip

We give more details for the case $n=0$. Thus $K=\cK=\exp(524)$, and we would
like to show, using proposition~\ref{prop:sevenK} that all odd integers $9K/10 \le N \le K$
are sums of seven non-negative cubes. 
The routine described in the proof of Lemma~\ref{lem:kappadef} gives the following suitable
value for $\kappa$:
\begin{multline*}
\kappa=
3
\times
173
\times 179
\times 191
\times 197
\times 227
\times 233
\times 239
\times 251
\times 257
\times 263
\times 269
\times 281
\times 293 \\
\times 311
\times 317 
\times 347
\times 353
\times 359
\times 383
\times 389
\times 401
\times 419
\times 431
\times 443
\times 207443 
\, .
\end{multline*}
Table~\ref{table:bigsift} gives some of the details for the computation.
We take $\cW=\cW_{48}$. Then $\#\cW=\#\cW_0+\sum \# \cU_{i}=9943$, and
\begin{align*}
\ell(\fJ(\cW) \cap [0,M))= &
12459371373955496388240157141404031514014113708989680551759\\
& 37887691670913319978\frac{1}{2}
\approx 1.25 \times 10^{78}.
\end{align*}
In comparison,
\[
M=M_{48} \approx 1.64 \times 10^{235}, \qquad K=3.72\times 10^{227},
\qquad  
\]
\begin{table}
\caption{This table gives details for the computation for the case $n=0$. For each $i \ge 1$,
our script computes $\fJ(\cW_i)$ as a disjoint union of subintervals of $[0,M_i)$. The number
of intervals is given in the fourth column. The fifth column gives, to $3$ significant figures,
the ratio $\ell_i/M_i$ where $\ell_i=\ell(\fJ(\cW_i) \cap [0,M_i))$.
}. 
\label{table:bigsift}
\centering
\begin{minipage}{0.47\textwidth}\centering
{
\tabulinesep=1.2mm
\begin{tabu}{||c|c|c|c|c||}
\hline\hline
& & & \small{number} & \\
$i$ & $p_{i-1}$ & $\# \cU_{i-1}$ & \small{of} & $\ell_i/M_i$\\ 
& & & \small{intervals} & \\
\hline\hline
$0$ & 
--
&
--
&
$854$
&
$1.90 \times 10^{-20}$
\\
\hline
$ 1 $ &
$ 449 $ &
$ 174 $ &
$ 775 $ &
$ 3.73\times 10^{-23} $\\ 
\hline
$ 2 $ &
$ 461 $ &
$ 175 $ &
$ 745 $ &
$ 7.94\times 10^{-26} $ \\
\hline
$ 3 $ &
$ 467 $ &
$ 176 $ &
$ 740 $ &
$ 1.70\times 10^{-28} $ \\
\hline
$ 4 $ &
$ 479 $ &
$ 176 $ &
$ 735 $ &
$ 3.54\times 10^{-31} $ \\
\hline
$ 5 $ &
$ 491 $ &
$ 177 $ &
$ 732 $ &
$ 7.20\times 10^{-34} $ \\
\hline
$ 6 $ &
$ 503 $ &
$ 178 $ &
$ 730 $ &
$ 1.42\times 10^{-36} $ \\
\hline
$ 7 $ &
$ 509 $ &
$ 178 $ &
$ 730 $ &
$ 2.80\times 10^{-39} $ \\
\hline
$ 8 $ &
$ 521 $ &
$ 179 $ &
$ 730 $ &
$ 5.38\times 10^{-42} $ \\
\hline
$ 9 $ &
$ 557 $ &
$ 181 $ &
$ 730 $ &
$ 9.65\times 10^{-45} $ \\
\hline
$ 10 $ &
$ 563 $ &
$ 181 $ &
$ 731 $ &
$ 1.71\times 10^{-47} $\\ 
\hline
$ 11 $ &
$ 569 $ &
$ 181 $ &
$ 730 $ &
$ 3.01\times 10^{-50} $\\ 
\hline
$ 12 $ &
$ 587 $ &
$ 182 $ &
$ 729 $ &
$ 5.13\times 10^{-53} $ \\
\hline
$ 13 $ &
$ 593 $ &
$ 182 $ &
$ 729 $ &
$ 8.64\times 10^{-56} $ \\
\hline
$ 14 $ &
$ 599 $ &
$ 183 $ &
$ 729 $ &
$ 1.44\times 10^{-58} $ \\
\hline
$ 15 $ &
$ 617 $ &
$ 183 $ &
$ 729 $ &
$ 2.34\times 10^{-61} $ \\
\hline
$ 16 $ &
$ 641 $ &
$ 185 $ &
$ 729 $ &
$ 3.64\times 10^{-64} $ \\
\hline
$ 17 $ &
$ 647 $ &
$ 185 $ &
$ 729 $ &
$ 5.63\times 10^{-67} $ \\
\hline
$ 18 $ &
$ 653 $ &
$ 185 $ &
$ 729 $ &
$ 8.62\times 10^{-70} $ \\
\hline
$ 19 $ &
$ 659 $ &
$ 185 $ &
$ 729 $ &
$ 1.31\times 10^{-72} $ \\
\hline
$ 20 $ &
$ 677 $ &
$ 186 $ &
$ 729 $ &
$ 1.93\times 10^{-75} $ \\
\hline
$ 21 $ &
$ 683 $ &
$ 186 $ &
$ 729 $ &
$ 2.83\times 10^{-78} $ \\
\hline
$ 22 $ &
$ 701 $ &
$ 187 $ &
$ 729 $ &
$ 4.04\times 10^{-81} $ \\
\hline
$ 23 $ &
$ 719 $ &
$ 188 $ &
$ 729 $ &
$ 5.61\times 10^{-84} $ \\
\hline
$ 24 $ &
$ 743 $ &
$ 189 $ &
$ 729 $ &
$ 7.55\times 10^{-87} $ \\
\hline
\end{tabu}
}
\end{minipage}
\hfill
\begin{minipage}{0.47\textwidth}\centering
{
\tabulinesep=1.2mm
\begin{tabu}{||c|c|c|c|c||}
\hline\hline
& & & \small{number} & \\
$i$ & $q_{i-1}$ & $\# \cU_{i-1}$ & \small{of} & $\ell_i/M_i$\\ 
& & & \small{intervals} & \\
\hline\hline
$ 25 $ &
$ 761 $ &
$ 189 $ &
$ 729 $ &
$ 9.93\times 10^{-90} $ \\
\hline
$ 26 $ &
$ 773 $ &
$ 190 $ &
$ 729 $ &
$ 1.28\times 10^{-92} $ \\
\hline
$ 27 $ &
$ 797 $ &
$ 191 $ &
$ 729 $ &
$ 1.61\times 10^{-95} $ \\
\hline
$ 28 $ &
$ 809 $ &
$ 191 $ &
$ 729 $ &
$ 1.99\times 10^{-98} $ \\
\hline
$ 29 $ &
$ 821 $ &
$ 192 $ &
$ 729 $ &
$ 2.43\times 10^{-101} $ \\
\hline
$ 30 $ &
$ 827 $ &
$ 192 $ &
$ 729 $ &
$ 2.93\times 10^{-104} $ \\
\hline
$ 31 $ &
$ 839 $ &
$ 192 $ &
$ 729 $ &
$ 3.50\times 10^{-107} $ \\
\hline
$ 32 $ &
$ 857 $ &
$ 193 $ &
$ 729 $ &
$ 4.08\times 10^{-110} $ \\
\hline
$ 33 $ &
$ 863 $ &
$ 193 $ &
$ 729 $ &
$ 4.73\times 10^{-113} $ \\
\hline
$ 34 $ &
$ 881 $ &
$ 194 $ &
$ 729 $ &
$ 5.36\times 10^{-116} $\\ 
\hline
$ 35 $ &
$ 887 $ &
$ 194 $ &
$ 729 $ &
$ 6.05\times 10^{-119} $\\
\hline
$ 36 $ &
$ 911 $ &
$ 195 $ &
$ 729 $ &
$ 6.64\times 10^{-122} $ \\
\hline
$ 37 $ &
$ 929 $ &
$ 195 $ &
$ 729 $ &
$ 7.14\times 10^{-125} $\\
\hline
$ 38 $ &
$ 941 $ &
$ 195 $ &
$ 729 $ &
$ 7.59\times 10^{-128} $\\
\hline
$ 39 $ &
$ 947 $ &
$ 196 $ &
$ 729 $ &
$ 8.02\times 10^{-131} $\\
\hline
$ 40 $ &
$ 953 $ &
$ 196 $ &
$ 729 $ &
$ 8.41\times 10^{-134} $ \\
\hline
$ 41 $ &
$ 971 $ &
$ 196 $ &
$ 729 $ &
$ 8.66\times 10^{-137} $\\
\hline
$ 42 $ &
$ 977 $ &
$ 197 $ &
$ 729 $ &
$ 8.87\times 10^{-140} $ \\
\hline
$ 43 $ &
$ 983 $ &
$ 197 $ &
$ 729 $ &
$ 9.02\times 10^{-143} $\\
\hline
$ 44 $ &
$ 1013 $ &
$ 198 $ &
$ 729 $ &
$ 8.91\times 10^{-146} $\\
\hline
$ 45 $ &
$ 1019 $ &
$ 198 $ &
$ 729 $ &
$ 8.74\times 10^{-149} $\\
\hline
$ 46 $ &
$ 1031 $ &
$ 198 $ &
$ 729 $ &
$ 8.48\times 10^{-152} $\\
\hline
$ 47 $ &
$ 1049 $ &
$ 199 $ &
$ 729 $ &
$ 8.08\times 10^{-155} $\\
\hline
$ 48 $ &
$ 1061 $ &
$ 199 $ &
$ 729 $ &
$ 7.62\times 10^{-158} $\\
\hline
\end{tabu}
}
\end{minipage}
\end{table}
The set $\fS$ as in \eqref{eqn:decomp}
turns out be precisely the set of $171$ rationals
$a/q \in [0,1]$ with denominators $q$ belonging to
\[
1,\; 2,\; 3,\; 4,\; 5,\; 6,\; 7,\; 8,\; 9,\; 10,\; 12,\;
 13,\; 14,\; 15,\; 16,\; 18,\; 19,\; 21,\; 24,\;
 26,\; 28,\; 30,\; 36,\; 42 \, .
\] 
As a check on our results, we apply Proposition~\ref{prop:ripple2} to show that there
is an interval close to $(a/42)\cdot M$ for $1 \le a \le 41$ with $\gcd(a,42)=1$.
Our $\cW$ and $M$ satisfy the hypotheses of Section~\ref{sec:ripples} with 
$L=K^{1/3}$.
Note that $3 \mid m \mid M$ for all $m \in \cW$. As $M$ is squarefree, we have $3 \nmid (M/m)$.
Moreover, all the prime divisors of $M/3$ are $\equiv 5 \pmod{6}$. It follows that $\gcd(aM/m,42)=1$
for all $m \in \cW$. Let $q=42$ and $s=d=5$ in Proposition~\ref{prop:ripple2}; hypothesis~\eqref{eqn:rineq}
is trivially satisfied. Now
$s+1,\dotsc,s+d$ are the integers $6$, $7$, $8$, $9$, $10$ and none of these 
are coprime to $42$. Thus condition~\eqref{eqn:aMm} is also satisfied. By Proposition~\ref{prop:ripple2},
for each $1 \le a \le 41$ with $\gcd(a,42)=1$ we have
\begin{equation}\label{eqn:42}
\left( \frac{a}{42}M -\frac{4471}{10500} \cdot K^{1/3}\,
, \, \frac{a}{42}M -\frac{73}{210} \cdot K^{1/3} \right)
\; \subseteq  \;
\fJ(\cW). 
\end{equation}
One of the $729$ intervals that make up $\fJ(\cW)$ is $[u,v)$ where the end points $u$, $v$
are
\begin{align*}
u= & 3895173640423584874713349032421520960246664873653293975537400307\\
   & 5224581570150578968661382487115397667257923729694373737120676906\\
   & 3932017310777324617938079775100516093231041460322490961793995991\\
   & 410145937421686204642056677472293123392066\frac{3}{10} \,  ,\\
v= & 3895173640423584874713349032421520960246664873653293975537400307\\
   & 5224581570150578968661382487115397667257923729694373737120676906\\
   & 3932017310777324617938079775101078697615391607190077739469387928\\
   & 238665618669912989320140106379011502569660 \, , \\
\end{align*}
and we checked that the interval in \eqref{eqn:42} with $a=1$ is contained in $[u,v)$.
It is also interesting to note how close the two intervals are in length: the 
ratio of the lengths of the two intervals is
\[
\frac{\left( 4471/10500 - 73/210 \right) \cdot K^{1/3}}{v-u} \approx 0.9994
\]
which illustrates how remarkably accurate our Proposition~\ref{prop:ripple2} is.

\section{Completing the Proof of Theorem~\ref{thm:main}}\label{sec:proofsmall}
It remains to apply Proposition~\ref{prop:sevenK}
to the intervals $(9/10)^{n+1} \cK \le N \le (9/10)^n \cK$
with $3998 \le n \le 4226$. 
We write $K=K_2=(9/10)^n \cK$ and $K_1=(9/10)^{n+1} \cK$.
It is no longer practical
to use the choices in Lemma~\ref{lem:parachoice} as the
interval in \eqref{eqn:parainterval} is too short
to contain many squarefree $m$ whose prime divisors are $3$
and small primes $\equiv 5 \pmod{6}$. The interval in \eqref{eqn:parainterval}
is a result of imposing the uniform choices $\varepsilon_m=0$ and $\delta_m=1/10$.
Instead we consider integers $m$ satisfying conditions (i)--(iii) of Section~\ref{sec:criterion}
but belonging to the (much larger) interval
\begin{equation}\label{eqn:secondrange}
\frac{12}{5}  K^{1/3} \le m \le \frac{16}{5} K^{1/3}\, .
\end{equation}
For each such $m$ we take 
$\varepsilon_m=\varepsilon^{\prime}/1000$
and
$\delta_m=\delta^\prime/1000$ where $\varepsilon^\prime$, $\delta^\prime$
are integers with $\varepsilon^\prime$ and $\delta^\prime$
respectively as small and as large as possible such that the 
conditions (v), (vi) of Section~\ref{sec:criterion} are satisfied.
We only keep those values of $m$ for which
\begin{equation}\label{eqn:restriction}
0 \le \varepsilon_m < \delta_m \le 1, \qquad \delta_m-\varepsilon_m \ge 1/20;
\end{equation}
an elementary though lengthy analysis in fact shows that the inequalities in \eqref{eqn:restriction}
together with (v) and (vi) force $m$ to belong to the interval \eqref{eqn:secondrange}.
Note that the set $\fJ(m,\varepsilon_m,\delta_m)$ has \lq relative density\rq\ $1-\delta_m+\varepsilon_m$
in $\R$; the restriction $\delta_m-\varepsilon_m \ge 1/20$ ensures that this relative density
is not too close to $1$, and that therefore $m$ makes a significant contribution to depleting
the intervals in Algorithms \ref{alg1} and \ref{alg2}.

We choose a prime $q \equiv 5 \pmod{6}$, 
depending on $K$, and let
\[
M_0=3 \cdot 5 \cdot 11 \cdot \cdots \cdot q
\]
which is the product of $3$ and the primes $\le q$ that are $\equiv 5 \pmod{6}$.
Let $\cW_0$ be the set of positive integers dividing $M_0$
and satisfying the above conditions. 
We found experimentally that for each $n$ in the above range
it is always possible to choose $q$ so that
\[
M_0=\lcm(\cW_0),
\qquad \prod_{m \in \cW_0} (1-\delta_m+\varepsilon_m) \le 1/5, \qquad 
\log_{10}(M_0/K^{1/3}) \le 7.5 \, .
\]
The inequality $\prod_{m \in \cW_0} (1-\delta_m+\varepsilon_m) \le 1/5$
indicates that $\ell(\fJ(\cW_0) \cap [0,M_0))$ should heuristically be at most $M_0/5$
which means that this is a good first step at depleting the interval
$[0,M_0)$. The other inequality indicates that we can 
compute $\fJ(\cW_0) \cap [0,M_0)$ in a reasonable number of steps according
to the heuristic following Algorithm~\ref{alg1}.
We let $p_0$ be the first prime $\equiv 5 \pmod{6}$ that is $>q$,
and $p_1$ the next such prime and so on. We let $M_{i+1}=p_{i} M_{i}$
and construct a tower as before. We stop once
$\fJ(\cW_i) \cap [0,M_i)$ satisfies the criterion of Proposition~\ref{prop:sevenK}.
Our \texttt{Magma} script succeeded in doing this for all $n$ in
the range $3998 \le n \le 4226$. 
The total CPU time was around $2750$ hours, 
but the computation was spread over $59$ processors so the 
actual time was less than $2$ days.

\bigskip

We give a few of details for the computation for the value $n=4226$.
The final $M$ is the product of $3$ and the primes $p \equiv 5 \pmod{6}$
that are $\le 227$. 
The final $\cW$ has $8083$ elements.
It turns out that $\fJ(\cW) \cap [0,M)$ consists of 
$305$ intervals and that
$\ell(\fJ(\cW) \cap [0,M))/M \approx 2.24 \times 10^{-32}$.

\bibliographystyle{plainnat}
\bibliography{cubesum}

\begin{thebibliography}{20}
\providecommand{\natexlab}[1]{#1}
\providecommand{\url}[1]{\texttt{#1}}
\expandafter\ifx\csname urlstyle\endcsname\relax
  \providecommand{\doi}[1]{doi: #1}\else
  \providecommand{\doi}{doi: \begingroup \urlstyle{rm}\Url}\fi

\bibitem[Baer(1913)]{Baer}
W.~S. Baer.
\newblock \"{U}ber die {Z}erlegung der ganzen {Z}ahlen in sieben {K}uben.
\newblock \emph{Math. Ann.}, 74\penalty0 (4):\penalty0 511--514, 1913.

\bibitem[Bertault et~al.(1999)Bertault, Ramar{\'e}, and Zimmermann]{BRZ}
F.~Bertault, O.~Ramar{\'e}, and P.~Zimmermann.
\newblock On sums of seven cubes.
\newblock \emph{Math. Comp.}, 68\penalty0 (227):\penalty0 1303--1310, 1999.

\bibitem[Boklan and Elkies(2009)]{BE}
K.~O. Boklan and N.~D. Elkies.
\newblock Every multiple of $4$ except $212$, $364$, $420$, and $428$ is the
  sum of seven cubes.
\newblock \emph{ArXiv e-prints}, March 2009.
\newblock arXiv:0903.4503.

\bibitem[Bosma et~al.(1997)Bosma, Cannon, and Playoust]{magma}
W.~Bosma, J.~Cannon, and C.~Playoust.
\newblock The {M}agma algebra system. {I}. {T}he user language.
\newblock \emph{J. Symbolic Comput.}, 24\penalty0 (3-4):\penalty0 235--265,
  1997.
\newblock Computational algebra and number theory (London, 1993).

\bibitem[Deshouillers et~al.(2000)Deshouillers, Hennecart, and Landreau]{Desh}
J.-M. Deshouillers, F.~Hennecart, and B.~Landreau.
\newblock {$7\,373\,170\,279\,850$}.
\newblock \emph{Math. Comp.}, 69\penalty0 (229):\penalty0 421--439, 2000.

\bibitem[Dickson(1927)]{DicksonExtension}
L.~E. Dickson.
\newblock Extensions of {W}aring's {T}heorem on {N}ine {C}ubes.
\newblock \emph{Amer. Math. Monthly}, 34\penalty0 (4):\penalty0 177--183, 1927.

\bibitem[Dickson(1939)]{DicksonEight}
L.~E. Dickson.
\newblock All integers except {$23$} and {$239$} are sums of eight cubes.
\newblock \emph{Bull. Amer. Math. Soc.}, 45:\penalty0 588--591, 1939.

\bibitem[Elkies(2010)]{Elkies}
N.~D. Elkies.
\newblock Every even number greater than 454 is the sum of seven cubes.
\newblock \emph{ArXiv e-prints}, September 2010.
\newblock arXiv:1009.3983.

\bibitem[Jacobi(1851)]{Jacobi}
C.~G.~J. Jacobi.
\newblock \"{U}ber die zusammensetzung der zahlen aus ganzen positiven cuben;
  nebst einer tabelle für die kleinste cubenanzahl, aus welcher jede zahl bis
  $12000$ zusammengesetzt werden kann.
\newblock \emph{Journal f\"{u}r die reine und angewandte Mathematik}, XLII,
  1851.

\bibitem[Kempner(1912)]{Kempner}
A.~Kempner.
\newblock Bemerkungen zum {W}aringschen {P}roblem.
\newblock \emph{Math. Ann.}, 72\penalty0 (3):\penalty0 387--399, 1912.

\bibitem[Landau(1911)]{Landau}
E.~Landau.
\newblock {\"Uber die Zerlegung positiver ganzer Zahlen in positive Kuben.}
\newblock \emph{{Arch. der Math. u. Phys. (3)}}, 18:\penalty0 248--252, 1911.

\bibitem[Linnik(1943)]{Linnik}
Yu.~V. Linnik.
\newblock On the representation of large numbers as sums of seven cubes.
\newblock \emph{Rec. Math. [Mat. Sbornik] N. S.}, 12(54):\penalty0 218--224,
  1943.

\bibitem[Maillet(1895)]{Maillet}
E.~Maillet.
\newblock Sur la decomposition d'un nombre entier en une somme de cubes
  d'entiers positifs.
\newblock \emph{Assoc. Fran\c{c}. Bordeaux}, XXIV:\penalty0 242--247, 1895.

\bibitem[McCurley(1984)]{McCurley}
K.~S. McCurley.
\newblock An effective seven cube theorem.
\newblock \emph{J. Number Theory}, 19\penalty0 (2):\penalty0 176--183, 1984.

\bibitem[Ramar{\'e}(2005)]{RamareBig}
O.~Ramar{\'e}.
\newblock An explicit seven cube theorem.
\newblock \emph{Acta Arith.}, 118\penalty0 (4):\penalty0 375--382, 2005.

\bibitem[Ramar{\'e}(2007)]{RamareSmall}
O.~Ramar{\'e}.
\newblock An explicit result of the sum of seven cubes.
\newblock \emph{Manuscripta Math.}, 124\penalty0 (1):\penalty0 59--75, 2007.

\bibitem[Romani(1982)]{Romani}
F.~Romani.
\newblock Computations concerning {W}aring's problem for cubes.
\newblock \emph{Calcolo}, 19\penalty0 (4):\penalty0 415--431, 1982.

\bibitem[von Sterneck(1903)]{Sterneck}
R.~D. von Sterneck.
\newblock {\"Uber die kleinste Anzahl Kuben, aus welchen jede Zahl bis 40000
  zusammengesetzt werden kann.}
\newblock \emph{{Wien. Ber.}}, 112:\penalty0 1627--1666, 1903.

\bibitem[Watson(1951)]{Watson}
G.~L. Watson.
\newblock A proof of the seven cube theorem.
\newblock \emph{J. London Math. Soc.}, 26:\penalty0 153--156, 1951.

\bibitem[Wieferich(1908)]{Wieferich}
A.~Wieferich.
\newblock Beweis des {S}atzes, da\ss\ sich eine jede ganze {Z}ahl als {S}umme
  von h\"ochstens neun positiven {K}uben darstellen l\"a\ss t.
\newblock \emph{Math. Ann.}, 66\penalty0 (1):\penalty0 95--101, 1908.

\end{thebibliography}

\end{document}